\documentclass[trsc,nonblindrev]{informs3} 

\OneAndAHalfSpacedXI 

\usepackage[colorlinks,
            linkcolor=blue,
            anchorcolor=blue,
            citecolor=blue
            ]{hyperref}
\usepackage{adjustbox}
\usepackage{booktabs}
\usepackage{subcaption}
\usepackage{xcolor}

\usepackage{natbib}
 \def\bibfont{\small}%
\TheoremsNumberedThrough     

\EquationsNumberedThrough    

\begin{document}




\TITLE{Dynamic Pricing for Air Cargo Revenue Management}

\ARTICLEAUTHORS{%
\AUTHOR{Chengyu Du}
\AFF{Department of Industrial Engineering, Tsinghua University}
\AUTHOR{Fang He\thanks{Corresponding author: fanghe@tsinghua.edu.cn}}
\AFF{Department of Industrial Engineering, Tsinghua University}
\AUTHOR{Xi Lin}
\AFF{Department of Civil and Environmental Engineering, University of Michigan, Ann Arbor}
} 

\ABSTRACT{%
We address a dynamic pricing problem for airlines aiming to maximize expected revenue from selling cargo space on a single-leg flight. The cargo shipments' weight and volume are uncertain and their precise values remain unavailable at the booking time. We model this problem as a Markov decision process, and further derive a necessary condition for its optimal pricing strategy. To break the curse of dimensionality, we develop two categories of approximation methods and pricing strategies. One category is based on the quantity of accepted bookings, while the other is founded on the expected weight and volume of accepted bookings. We prove that the pricing strategy of the quantity-based method possesses several inherent structural properties, which are crucial for analytically validating the model and accelerating the computational process. For the weight-volume-based approximation method, we derive a theoretical upper bound for the optimality gap of total expected revenue. For both methods, we further develop augmented strategies to address the extreme pricing issues in scenarios with high product heterogeneity and incorporate the second moment to enhance performance in the scenarios of high uncertainty, respectively. We utilize realistic dataset to conduct extensive numerical tests, and the results show that the average performance gap between the optimal expected revenue and that of each proposed pricing strategy is less than 10\%. The quantity-based method requires the least computation, and performs quite well in the scenarios with low product heterogeneity. The augmented quantity-based method further promotes the performance resilience to product heterogeneity. The weight-volume-based method further enhances the resilience to product heterogeneity, and exhibits superior performance in the scenarios with low uncertainty. The augmented weight-volume-based method significantly improves the revenue when there are high penalties for overbooking and high uncertainty.
}%



\maketitle

%


\section{Introduction}
\label{sec:intro}
Air cargo plays a vital role in contemporary logistics and serves as a critical infrastructure for facilitating global trade. Despite constituting just 1\% of global shipments in terms of volume, it contributes to approximately 35\% of the total value \citep{IATA2015}. In recent years, driven by the rapid expansion of the aviation sector and the maturation of air transportation networks, the air cargo industry has achieved an average annual growth rate of 3.6\%. Projections suggest that this growth will persist at a rate of 4.1\% annually for the next two decades \citep{Boeing2022}.

Shippers, freight forwarders, and airlines (carriers) are three principal air cargo logistics network participants.
Freight forwarders should secure cargo capacity from airlines to meet shippers' demands.
Airlines allocate cargo capacity through a two-phase process. In the initial phase, which takes place six months to a year before flight departure, freight forwarders submit bids to airlines to secure cargo capacity. Airlines then decide which allotment contracts to sign. In the second phase, occurring approximately two weeks before flight departure, freight forwarders must confirm their allocated capacity. Any remaining capacity becomes available for sale by the airlines in the spot market \citep{Popescu2006}. 

Revenue management (RM), which originated in the air passenger industry, has experienced rapid development as a means to optimize resource utilization and maximize revenue since the passing of the Airline Deregulation Act in 1978 in the United States.
Estimates suggest that RM has led to a revenue increase of approximately 4-5\% for airlines \citep{Talluri2004}. 
While RM has been successfully implemented by most airlines in the air passenger industry, its application in the air cargo industry continues to present challenges. There are some principal distinctions between RM of air cargo and air passenger  \citep{Amaruchkul2007, barz2016air}. One of the primary differences lies in dimensionality. Air passenger RM primarily deals with one-dimensional capacity (seat number), whereas air cargo involves multi-dimensional capacity, including weight and volume. Moreover, owing to the flexibility of cargo booking requests, the precise weight and volume remain indeterminate before flight departure. In other words, at the booking time, the weight and volume of cargo remain uncertain.

Among the studies focused on RM in industrial applications, such as air transportation, there are two approaches of implementation, capacity control and dynamic pricing. Capacity control refers to managing bookings by accepting or rejecting them based on the presumption that demand is independent \citep{klein2020review}. The majority of existing research on air cargo RM adopts capacity control \citep[e.g.][]{pak2004cargo, Amaruchkul2007, Huang2010, hoffmann2013dynamic}. Among them, \cite{pak2004cargo} used bid-price control (BPC), which is a subcategory of the capacity control. BPC involves setting bid (threshold) prices for resources or inventory units, such as weight and volume capacities or airline seats. A booking is only accepted if the offered fare surpasses the summation of bid prices for all required resources \citep{talluri1998analysis}. 
Dynamic pricing is another fundamental and widely-employed approach in RM. Under the assumption of independent demand, the seller determines dynamic prices of products over time and available resources, while the customers make purchase decisions based on these prices. Different from the capacity control, in dynamic pricing, the price serves as the principal decision variable. Throughout the sales horizon, customers' maximum willingness to pay, commonly referred to as their "reservation prices", may fluctuate due to the evolving urgency of their booking needs. Notably, dynamic pricing offers enhanced adaptability to these temporal shifts in consumer valuation, potentially leading to augmented revenue streams. Despite its merit, the studies on the application of dynamic pricing in cargo RM are rather limited. The possible reason is that the multi-dimensional characteristics and uncertainty of weight and volume of cargo bookings make it challenging to address the curse of dimensionality and develop robust pricing strategies to maximize revenue. To fill this gap, this research focuses on the implementation of dynamic pricing in air cargo RM for a single-leg flight during the spot sale phase to maximize revenue for the airline. Our contributions can be summarized as follows.

\begin{itemize}
    \item We develop a Markov decision process (MDP) to model
    the dynamic pricing problem during air cargo spot sales for a single-leg flight with multiple booking types, and further derive a necessary condition for its optimal pricing strategy. The type of a booking request is represented through the weight and volume distributions, as well as the distribution of reservation prices. The precise weight and volume of a booking become determinant only at the end of the booking horizon. Reservation price indicates customers’ highest willingness to pay for the booking. And bookings with the same weight-volume distribution can have different distributions of reservation prices due to attributes such as the perishability of the goods and packaging quality. 
    \item We develop two types of quantity-based approximation pricing methods to break the curse of dimensionality. The main idea is aggregating the multi-dimensional state space into a one-dimensional space. We prove that the pricing strategy of the primal method satisfies several structure properties, which are of paramount importance for analytically validating the proposed model, as the consistency between analytical and empirical observations serves as a robust indicator of the model's validity. Moreover, the monotonicity of the pricing strategy can also reduce the searching domain during the iterative process of computing the optimal strategy. The second type of method utilizes an augmented pricing strategy to further address the extreme pricing issues in scenarios with high product heterogeneity. 
    \item To further enhance the performance of dynamic pricing in high product-heterogeneity scenario, we propose two types of weight-volume-based approximation pricing methods. Specifically, we first formulate a new MDP characterized by states reflecting expected weight and volume. Subsequently, we discretize this MDP and derive its optimal value function through backward induction. The derived optimal value function is then utilized to compute a weight-volume-based approximation pricing strategy. Furthermore, we derive a theoretical bound concerning the optimality gap of the total expected revenue. To bolster performance under scenarios marked by pronounced uncertainty, the second type of method embraces an augmented pricing strategy, incorporating the second moments of weight and volume distributions associated with diverse booking types.
    \item We undertake comprehensive numerical tests utilizing a realistic dataset. The outcomes indicate that both quantity-based and weight-volume-based approximation methods yield superior performance, with the maximum average deviation from the optimal expected revenue remaining below 10\%. The quantity-based method requires less computation, and approaches optimal expected revenue when there is low heterogeneity in the weight and volume distributions of diverse booking types. Its augmented variant exhibits more robust performance with respect to the heterogeneity of booking types. While the weight-volume-based approach is computationally more intensive, it displays superior performance in scenarios characterized by pronounced heterogeneity of booking types. Furthermore, we show that the performance improvement from its augmented variant by incorporating second moments increases as the uncertainty regarding weight and volume or the penalty for overbooking increases.
\end{itemize}

The rest of this paper is organized as follows. In Section \ref{sec:review}, we review the relevant literature. In Section \ref{sec:pro}, we provide the model formulation. In Section \ref{sec:quan}, we develop quantity-based approximation methods and discuss their properties. Section \ref{sec:weight}  proposes weight-volume-based approximation methods and derives an upper bound of revenue gap for the primal weight-volume-based method. Section \ref{sec:NE} conducts numerical experiments. Finally, in Section \ref{sec:conclude}, we conclude.

\section{Literature Review} \label{sec:review}
The major research on air cargo operations explores flight scheduling, fleet routing, and RM problems. 
Refer to \cite{feng2015air} for an informative review of air cargo operations. intricacies
Air cargo flight scheduling and fleet routing comprise route choice, assignments of cargo to flights, and fleet and crew scheduling \citep{yan2006air,bartodziej2009models,azadian2012dynamic,lee2019design,xiao2022integrated,zheng2023air}. For instance, \cite{yan2006air} proposed a planning framework integrating airport choice, fleet routing and timetable development to optimize operating profit. \cite{azadian2012dynamic} formulated an MDP to route air cargo dynamically considering live data. \cite{lee2019design} proposed a comprehensive framework that addressed both the design of a cargo transportation network and the choice of fleet routes.
As our work directly relates to the operational decisions concerning RM, we primarily review literature in RM field.

Extensive literature has researched the domain of passenger RM \citep[e.g.][]{belobaba1987survey, mcgill1999revenue, Talluri2004}. However, the realm of cargo RM remains relatively under-researched and warrants further scholarly exploration. Several papers provide a practical perspective on air cargo RM. \cite{Kasilingam1997} compared passenger RM with cargo RM and introduced the attributes and intricacies of the latter. \cite{Billings2003} described the components of the cargo RM system from a business perspective and discussed the challenges faced in implementing RM. \cite{slager2004implementation} described how Koninklijke Luchtvaart Maatschappij (KLM) implemented air cargo RM and key factors for successful implementation. From a commercial standpoint, \cite{Becker2007} offered an enhanced list of complexities based on the supply side and the demand side for air cargo RM.

A number of studies delve into the cargo RM with a static formulation \citep[e.g.][]{karaesmen2001three, kasilingam1997economic, luo2009two}. Other papers model the booking process as a dynamic formulation with capacity control. 
\cite{Amaruchkul2007} demonstrated in this regard with their formulation of spot sales for single-leg cargo RM as a multi-dimensional MDP.
Their model assumes that the number of each type of booking is known in each period, while the volume and weight of each accepted booking remain indeterminate until just before the flight departs. To address overbooking, they introduced penalties for cargo unloading at the booking horizon's closure. They formulated heuristic policies and established bounds by decoupling the issue into one-dimensional subproblems of volume and weight, subsequently undertaking numerical analyses. Building on this foundation, in order to improve the heuristics presented in \cite{Amaruchkul2007}, a heuristic based on sampling approximation is developed by \cite{Huang2010}. 
In addition, \cite{hoffmann2013dynamic} developed a heuristic that offers monotone opportunity cost, which reduces computations and makes the manager can control capacity based on monotone switching curves. 
\cite{xiao2010revenue} proposed a continuous-time stochastic capacity control problem to model a two-dimensional RM issue. They derived an analytical solution and demonstrated that the optimal control strategy is threshold-based under specific conditions. A few papers use bid-price control in the cargo booking process \citep{pak2004cargo, han2010markov}. 
In terms of air cargo network RM, to optimize the capacity control policy, \cite{levina2011network} proposed a linear programming method based on simulation considering uncertain capacity.  
For a similar problem formulation, \cite{barz2016air} proposed several heuristic policies and engaged in the comparative analysis of the upper bounds derived from them. 
\cite{huang2015linear} developed deterministic and probabilistic linear programming models to estimate the bookings' opportunity cost. 
\cite{delgado2019multistage} approached the network capacity allocation as a multistage stochastic programming model for revenue maximization. 
In addition, \cite{previgliano2022managing} considered the capacity control network RM with random availability of physical capacities and time flexibility. They proposed a stochastic gradient algorithm to approximately solve it. 
Essentially, several features distinguish our work from the papers above. Primarily, our focus is on the dynamic pricing of different types of bookings in the spot market, whereas the papers above utilize a capacity control approach. Second, one of our proposed approximation pricing methods leverages the second-order information of uncertainty, a dimension not addressed by the heuristics in the previously mentioned papers.

Dynamic pricing differs from capacity control in that it allows for adjusting product prices over time.
\cite{chen2015recent} presented a comprehensive review on dynamic pricing problems in RM. 
The majority of research on pricing has concentrated on single-product scenarios. 
\cite{Gallego1994} considered the single-product dynamic pricing and showed that a constant price is near-optimal.
\cite{bitran1997periodic} addressed the problem of single-product pricing while considering customers with different reservation prices for the product.
\cite{zhao2000optimal} examined the problem in \cite{bitran1997periodic} and derived several structural properties.
\cite{feng1999maximizing} optimized the pricing of single commodity considering that some companies may be risk-sensitive. 
Several papers explore dynamic pricing involving multiple products or network structure.  
A multi-product dynamic pricing issue is investigated by \cite{gallego1997multiproduct}. 
They established a deterministic optimization problem and showed that the pricing strategy derived from it is near-optimal. 
A similar issue under network RM structure is studied by \cite{erdelyi2011using}. They proposed heuristics based on dynamic programming decomposition and showed that their approaches achieve great total expected revenue.
\cite{koenig2010list} conducted a novel comparison study between the use of capacity control and dynamic pricing with multiple products.
\cite{zhang2013assessing} conducted a comprehensive numerical experiment to assess the effect of dynamic pricing. They also adopted a decomposition method to compute the dynamic pricing policy, demonstrating its potential for revenue improvement. Additionally, for network pricing RM, \cite{ke2019approximate} proposed an approximate dynamic programming technique to optimize it. 
Our research differs from the papers above in that we consider multiple bookings with continuous and uncertain weight and volume in the air cargo industry. These features pose significant challenges in determining an effective solution, necessitating a pricing strategy that is sufficiently resilient to these inherent uncertainties.

\section{Problem Statement}
\label{sec:pro}

\subsection{Preliminaries}
We consider a single flight segment during spot sales with weight capacity $C_w$ and volume capacity $C_v$. The length of the booking horizon is $L$. Following \cite{Amaruchkul2007}, as well as \cite{levina2011network} and \cite{huang2015linear}, each booking is assumed to be categorized into one of $m$ distinct types, with the arrival being independent over time. Each booking type is characterized by its weight and volume distribution, packaging type, cargo type, and priority, which determines the customers' valuation of charges.
The arrival of type $i$ booking follows a Poisson process with a rate parameter of $\lambda_s^i$ ($s\in[0,L]$). Therefore, the arrival rate for all booking requests is $\lambda_s=\sum_{i=1}^m\lambda_{s}^i$. We segment the time horizon into $T$ periods, each of duration length $\Delta t$, so that there is at most one booking request arrives per period. $L=T\Delta t$, and we count the periods in the booking horizon forwards. That is, the booking horizon begins at period $t=0$ and ends with the flight's departure at period $t=T$. Type $i$ booking request corresponds to a joint distribution of weight $W_i$ and volume $V_i$. 
Upon the arrival of each booking request, the airline is aware of its type and corresponding weight and volume joint distribution. However, at the time of pricing, the airline is not aware of its precise weight and volume. When the flight is about to depart and the bookings are delivered, the airline becomes aware of the precise volume and weight.

In accordance with the International Air Transportation Association (IATA) guidelines, the charge for a type $i$ cargo with weight $w$ and volume $v$, is determined as $r_i\max\{w,v/\gamma\}$, where $\gamma=6000 \ \text{cm}^3/\text{kg}$; 
$\max\{w,v/\gamma\}$ is designated as the chargeable weight; and $r_i$ denotes the pricing for the unit chargeable weight of the cargo. 
Similar to \cite{bitran1997periodic}, it is assumed that each customer possesses a reservation price for unit costs in chargeable weight for each booking type. This represents customers’ highest willingness to pay for a booking. Should one customer's reservation price for her booking meet or exceed the pricing established by the airline, she will accept the booking at the current pricing. Otherwise, the airline will lose that booking. Assume that for each booking type $i$, a customer's reservation price is a stochastic variable, denoted as $RP_s^i$. Let $F_s^i(r)$ denote the corresponding cumulative distribution function, assumed to be differentiable, and $f_s^i(r)$ denote the corresponding probability distribution function at time $s\ (s\in[0,L])$. Distinct booking types may exhibit varied reservation price distributions. 
Upon the arrival of a booking of type $i$, it is assumed that the airline is only aware of the corresponding distribution of customers' reservation prices rather than a precise value. If the airline determines the price to be $r$ at period $t$, then the corresponding probability of acceptance for that booking is $1-F^i_{t\Delta t}(r)$. Note that it is assumed that reservation price distributions change smoothly across time, so we utilize $F^i_{t\Delta t}(r)$ to approximate it at period $t$. For convenience, we utilize the notation $F^i_{t}(r)$ to represent the cumulative distribution of reservation prices at period $t$.

We aim to optimize the revenue generated through spot sales by implementing dynamic pricing strategies for each type of booking.

\subsection{Model Formulation}
We model the spot sale dynamic pricing problem through an MDP and refer to it as the general model. We refer to $\mathbb{N}^m$ as the m-fold cross product of the non-negative integers. And we represent a vector in boldface.
The system's state $\mathbf{x}=(x_1, x_2,\dots, x_m)^T \in \mathbb{N}^m$, indicates the number of various accepted bookings, where $x_i$ denotes the quantity of accepted type $i$ bookings specifically.
At period $t$, we determine the pricing for type $i$ bookings based on state $\mathbf{x}$ to maximize expected revenue. We denote this pricing as $r_i(t,\mathbf{x})$, and for convenience, we refer to it as $r_i$.

A pricing strategy $\pi$ consists of $\mathbf{r}^{\pi}(t,\mathbf{x})=(r_1^{\pi}(t,\mathbf{x}), r_2^{\pi}(t,\mathbf{x}), \dots, r_m^{\pi}(t,\mathbf{x}))^T$ for all $\mathbf{x},\ t=0,1,\dots,T-1$. Let $\Pi$ be the set containing all pricing strategies. We define $J_t^{\pi}(\mathbf{x})$ as the expected revenue when the airline begins at state $\mathbf{x}$ at the start of period $t$ using pricing strategy $\pi$. $V_t(\mathbf{x})=\sup_{\pi\in\Pi}J_t^{\pi}(\mathbf{x})$ represents the maximum expected revenue. Note that we do not consider discounts because the bookings occur over a short period. The Bellman equation for the MDP model is given as Equation \eqref{eq:general}.

\begin{align}
    V_t(\mathbf{x})=&{\rm sup_{r_i}}E\Bigg\{\sum_{i=1}^mm_{t}^i(1-F_t^i(r_i))\left(r_i\max\left\{W_i, \frac{V_i}{\gamma}\right\}+V_{t+1}(\mathbf{x}+\mathbf{e}_i)\right) + \nonumber \\ &
    \quad\quad\quad \left(1-\sum_{i=1}^mm_{t}^i(1-F_t^i(r_i))\right)V_{t+1}(\mathbf{x}) \Bigg\} \nonumber\\
    =&{\rm sup_{r_i}}\Bigg\{\sum_{i=1}^mm_{t}^i(1-F_t^i(r_i))\left(r_iE\left[\max\left\{W_i, \frac{V_i}{\gamma}\right\}\right]+V_{t+1}(\mathbf{x}+\mathbf{e}_i)\right) + \nonumber\\
    & \quad\quad\quad \left(1-\sum_{i=1}^mm_{t}^i(1-F_t^i(r_i))\right)V_{t+1}(\mathbf{x}) \Bigg\} \nonumber\\
    =&{\rm sup_{r_i}}\Bigg\{\sum_{i=1}^mb_t^i(r_i)\left(r_iQ_i-(V_{t+1}(\mathbf{x})-V_{t+1}(\mathbf{x}+\mathbf{e}_i))\right) + V_{t+1}(\mathbf{x}) \Bigg\} \nonumber\\
    =&\sum_{i=1}^m{\rm sup_{r_i}}\big\{b_t^i(r_i)\left(r_iQ_i-(V_{t+1}(\mathbf{x})-V_{t+1}(\mathbf{x}+\mathbf{e}_i))\right)\big\} + V_{t+1}(\mathbf{x}), \nonumber\\
    &\quad\quad\quad\quad\quad\quad t=0,1,2,...,T-1.
\label{eq:general}
\end{align}
where $m_t^i=\int_{t\Delta t}^{(t+1)\Delta t}\lambda_{s}^ids$ denotes the arrival probability of a type $i$ booking at period $t$; $b_t^i(r_i)=m_{t}^i(1-F_t^i(r_i))$;
$Q_i=E\left[\max\left\{W_i, \frac{V_i}{\gamma}\right\}\right]$ denotes expected chargeable weight of a type $i$ booking; $\mathbf{e}_i$ denotes the $i$th unit m-vector.

The last equation in Equation \eqref{eq:general} indicates that the multi-dimensional optimal pricing optimization for each period can be decomposed into $m$ separate single-variable optimization problems. This decomposition simplifies the determination of optimal pricing for each period. For example, Newton's method can be applied to resolve each single-variable optimization issue. 

If the airline faces situations of overbooking, where the total weight or volume of accepted bookings surpasses the capacity, it has to decide which cargo to unload before the flight departs. The penalties incurred from unloading can be affected by various factors, such as the customer type, the cargo type, and associated charges, which can make prompt decision-making challenging \citep{Amaruchkul2007}. Following \cite{Amaruchkul2007} and \cite{Huang2010}, we assume the overbooking penalty is a separable function of the total weight and volume of accepted bookings.
Suppose $(W_{ij},V_{ij})$ represents the weight and volume of the $j$th type $i$ booking, with the distribution of $(W_{ij},V_{ij})$ being independent and identical to the distribution of $(W_i,V_i)$ for all $i$ and $j$.
The boundary condition at $t=T$ (taking off) is described by Equation \eqref{eq:gen_bound}.

\begin{equation}
    V_T(\mathbf{x})=-E\left[h_w\left((\sum_{i=1}^{m}\sum_{j=1}^{x_i}W_{ij}-C_w)^+\right)+h_v\left((\sum_{i=1}^{m}\sum_{j=1}^{x_i}V_{ij}-C_v)^+\right)\right],
\label{eq:gen_bound}
\end{equation}
where $h_w$ and $h_v$ are non-decreasing convex functions, with $h_w(0)=h_v(0)=0$. 

It is easy to derive two properties on $V_t(\mathbf{x})$: 
\begin{enumerate}
    \item[(1)] $V_t(\mathbf{x})\geq V_{t+1}(\mathbf{x}), \quad \forall \mathbf{x}\geq0, \ t=0,1,...,T-1$;
    \item[(2)] $V_t(\mathbf{x})\geq V_{t}(\mathbf{x}+\mathbf{e}_i), \quad \forall \mathbf{x}\geq0, \ i=1,2,...,m, \ t=0,1,...,T$.
\end{enumerate}

The following theorem provides a necessary condition for the optimal pricing of the MDP.

\begin{theorem}
A necessary condition for the optimal price $r_i$ at period $t$ and state $\mathbf{x}$ is:

\begin{equation}
    r_i = \frac{1-F_t^i(r_i)}{f_t^i(r_i)}+\frac{V_{t+1}(\mathbf{x})-V_{t+1}(\mathbf{x}+\mathbf{e}_i)}{Q_i}.
    \label{eq:ropt_gen}
\end{equation}
If the function $(1-F_t^i(r_i))^2/f_t^i(r_i)$ is decreasing in $r_i$, then the Equation \eqref{eq:ropt_gen} has a unique solution, which is the optimal price.
\label{prop:V_x} 
\end{theorem}

Throughout this study, without specific notification, the proofs are provided in Appendix.

When the distribution of the reservation price for type $i$ bookings follows an exponential or Weibull distribution, the condition $(1-F_t^i(r_i))^2/f_t^i(r_i)$ decreasing in $r_i$ holds, and the Equation \eqref{eq:ropt_gen} becomes a sufficient condition for the optimal pricing.

Equation \eqref{eq:ropt_gen} demonstrates $r_i(t,\mathbf{x}) \geq (V_{t+1}(\mathbf{x})-V_{t+1}(\mathbf{x}+\mathbf{e}_i))/Q_i$. The right term of this inequality, $(V_{t+1}(\mathbf{x})-V_{t+1}(\mathbf{x}+\mathbf{e}_i))/Q_i$, represents a type $i$ booking's unit opportunity cost at period $t$ and state $\mathbf{x}$. This indicates that regardless of the state, the optimal pricing is always at least equal to the opportunity cost.

Based on the Equations \eqref{eq:general}-\eqref{eq:gen_bound}, we aim to obtain the maximum expected revenue, denoted as $V_0(\mathbf{0})$, and the corresponding optimal pricing strategy. To solve these equations, we first need to define a maximum acceptable quantity of bookings, referred to as $x_{max}$. For instance, $x_{max}=C_w/\min_i\{E[W_i]\}$. We only assess the state $\mathbf{x}$ where the sum of its components, denoted as $||\mathbf{x}||_1$, does not exceed $x_{max}$. Then, we can employ backward induction to solve these equations. This is reasonable because of the finite time, arrival rate, and capacity. However, the computational load becomes intractable for most practical problems because we have to solve at least $O((x_{max})^mT)$ optimization problems during the backward induction. This leads to the challenge known as the curse of dimensionality, especially when $m$ is large. Thus, we propose several approximation methods in the following sections to address this challenge.

\section{Quantity-based Approximation Method}
\label{sec:quan}
When dealing with the curse of dimensionality in solving  MDPs, a commonly used technique for approximation is to aggregate the state and decision spaces. The primary idea behind the quantity-based approximation method is aggregating the state space and reducing the $m$-dimensional MDP to a one-dimensional MDP. We introduce the primal quantity-based approximation method and its structural properties in Section \ref{subsec:PQA}. Then in Section \ref{subsec:AQA}, its augmented variant is introduced.

\subsection{Primal Quantity-based Approximation Method}
\label{subsec:PQA}
We aggregate the state vector $\mathbf{x}$ in the general model by summing all its elements, denoted as $||\mathbf{x}||_1$. We use the notation $x$ to represent this aggregated value, corresponding to the total number of accepted bookings.

Parameter normalization is an important consideration in the quantity-based approximation method. 
To catch the stochastic weight and volume of bookings, we incorporate a parameter normalization utilizing the first and second moments of their distributions.
When a booking arrives, denote $p_i=\int_0^L\lambda_s^ids/\int_0^L\lambda_sds$ as the approximate probability that it is of type $i$.
We define a weight random variable $W^s$ that follows a normal distribution with mean $w^s=\sum_{i=1}^m p_iw_i$ and variance $\sigma^2_{W^s}=\sum_{i=1}^mp_i[\sigma_{W_i}^2+(w_i-w^s)^2]$, where $w_i$ and $\sigma^2_{W_i}$ are the mean and variance of the weight, respectively, for a type $i$ booking. Meanwhile, we define a volume random variable $V^s$ that follows a normal distribution with mean $v^s=\sum_{i=1}^m p_iv_i$ and variance $\sigma^2_{V^s}=\sum_{i=1}^mp_i[\sigma_{V_i}^2+(v_i-v^s)^2]$, where $v_i$ and $\sigma^2_{V_i}$ are the mean and variance of the volume, respectively, for a type $i$ booking. $W^s$ and $V^s$ represent the estimated weight and volume distribution for any booking.

Then, we formulate a new quantity-based MDP (Q-MDP). We use $V_t^{PQ}(x)$ to denote its optimal value function, and it can be derived iteratively through the Bellman equation and boundary condition:

\begin{equation}
\begin{aligned}
    V_t^{PQ}(x)=&\sum_{i=1}^m{\rm sup_{r_i}}\left\{b_t^i(r_i)\left(r_iQ_i-(V_{t+1}^{PQ}(x)-V_{t+1}^{PQ}(x+1))\right)\right\} + V_{t+1}^{PQ}(x), \\
    &\quad\quad\quad\quad\quad\quad t=0,1,2,...,T-1.
\end{aligned}
\label{eq:PQ}
\end{equation}
\begin{equation}
    V_T^{PQ}(x)=-E\left[h_w\left((\sum_{l=1}^xW_l^s-C_w)^+\right)+h_v\left((\sum_{l=1}^xV_l^s-C_v)^+\right)\right],
\label{eq:PQ_bound}
\end{equation}
where the distribution of $(W_{l}^{s},V_{l}^{s})$ is independent and identical to the distribution of $(W^s,V^s)$ for all $l$.

The primal quantity-based (PQ) approximation method approximates $V_t(\mathbf{x})$ by $V_t^{PQ}(||\mathbf{x}||_1)$. Let $\mathbf{r}^{PQ}(t,x)=(r_1^{PQ}(t,x),r_2^{PQ}(t,x),...,r_m^{PQ}(t,x))^T$ represent the solutions obtained by solving the Equations \eqref{eq:PQ}-\eqref{eq:PQ_bound}. 
The approximation method is equivalent to pricing a type $i$ booking at state $\mathbf{x}$ and period $t$ with the price $r_i^{PQ}(t,||\mathbf{x}||_1)$. This pricing method is referred to as PQ pricing. 
If there are multiple solutions, the pricing is determined by selecting the highest value. Similar to Theorem \ref{prop:V_x}, the PQ pricing $r_i^{PQ}(t,x)$ satisfies the following corollary:

\begin{corollary}
    A necessary condition for the price $r_i$ to be the solution of the Equations \eqref{eq:PQ}-\eqref{eq:PQ_bound} at period $t$ and state $x$ is:

    \begin{equation}
        r_i = \frac{1-F_t^i(r_i)}{f_t^i(r_i)}+\frac{V_{t+1}^{PQ}(x)-V_{t+1}^{PQ}(x+1)}{Q_i}.
        \label{eq:ropt_Q}
    \end{equation}
    If the function $(1-F_t^i(r_i))^2/f_t^i(r_i)$ is decreasing in $r_i$, then the Equation \eqref{eq:ropt_Q} has a unique solution.
    \label{prop:VQ_x} 
\end{corollary}

It is important to determine if an MDP strategy exhibits any patterns, such as monotonicity, due to its facilitation to implementation, attractiveness to decision-makers, and acceleration of computation \citep{puterman2014markov}. As such, we explore several structural properties of $V_t^{PQ}(x)$ and $r_i^{PQ}(t,x)$. First, we focus on $V_t^{PQ}(x)$. In the Q-MDP, at state $x$ and period $t$, the difference $V_t^{PQ}(x)-V_t^{PQ}(x+1)$ quantifies the opportunity cost of accepting a booking. The following propositions indicate it increases with the number of accepted bookings $x$ when $t$ is fixed and decreases over time with a fixed $x$. The proofs of them rely on the following lemma.

\begin{lemma}
    Assume function $g(x)$ is non-decreasing, and convex. $X$ is a non-negative random variable with finite mean and variance.
    For the convolution of $n$ independent and identically distributed random variables with the same distribution as $X$, denoted as $S_n = \sum_{i=1}^n X_i$, the following inequality holds:
    \begin{equation*}
        E[g(S_{n+1})] - E[g(S_n)] \geq E[g(S_n)] - E[g(S_{n-1})].
        \label{eq:cvx_lem}
    \end{equation*}
\label{lem:cvx}
\end{lemma}

\begin{proposition}
    $V_t^{PQ}(x)-V_t^{PQ}(x+1)\geq V_t^{PQ}(x-1)-V_t^{PQ}(x), \ t=0,1,...,T.$
\label{prop:inventory}
\end{proposition}

\begin{proposition}
    $V_t^{PQ}(x)-V_t^{PQ}(x+1)\geq V_{t+1}^{PQ}(x)-V_{t+1}^{PQ}(x+1), \ t=0,1,...,T-1.$
\label{prop:time}
\end{proposition}

Propositions \ref{prop:inventory}-\ref{prop:time} illustrate several properties of $V_t^{PQ}(x)$. Subsequently, the properties of PQ pricing are presented as follows.

\begin{proposition}
    $\mathbf{r}^{PQ}(t,x)\leq \mathbf{r}^{PQ}(t,x+1), \ t=0,1,...,T-1. $
\label{prop:Qprice_x}
\end{proposition}

Proposition \ref{prop:Qprice_x} states that PQ pricing increases with the number of accepted bookings at any fixed time. With regard to the monotonicity property of PQ pricing  with respect to time, we show in the following proposition that, under the condition of homogeneous arrival and reservation prices, PQ pricing decreases over time. We also note that this condition is sufficient but not necessary. 

\begin{proposition}
    If the distribution of reservation prices for each type of booking does not change over time and that the arrival of each type of booking is time-homogeneous, i.e., $F_s^i(r)=F^i(r)$ and $\lambda_s^i=\lambda^i$ for all $s\in [0,L]$, the following relationship holds:
    \begin{equation*}
        \mathbf{r}^{PQ}(t+1,x)\leq \mathbf{r}^{PQ}(t,x), \ t=0,1,...,T-2.
    \end{equation*}
\label{prop:Qprice_t}
\end{proposition}

Intuitively, if reservation price distribution shifts towards higher values as time progresses, PQ pricing will also increase. In such cases, Proposition \ref{prop:Qprice_t} no longer holds.


The structure proprieties of PQ pricing are essential for confirming PQ model's analytical validity. This is crucial as consistency between analytical findings and observed empirical behavior stands as a strong testament to the model's validity. Moreover, they are also useful for numerical computations as they can reduce the searching domain during backward induction. 

\begin{remark}
    PQ pricing is equivalent to optimal pricing if the weight and volume distributions of diverse booking types are independent and follow the same normal distribution.
\end{remark}

\subsection{Augmented Quantity-based Approximation Method}
\label{subsec:AQA}
While PQ pricing exhibits desirable structural properties, which enhance its intuitiveness and can reduce the search domain during backward induction, its performance may be suboptimal in certain scenarios due to its pricing structure. Concretely, PQ pricing method approximates all $V_t(\mathbf{x})$ by $V_t^{PQ}(x)$, where $||\mathbf{x}||_1=x$. According to Equations \eqref{eq:ropt_gen} and \eqref{eq:ropt_Q}, the difference between the optimal pricing equation and the PQ pricing equation lies in the second term on the right side. PQ pricing employs a single value, $V_{t+1}^{PQ}(x)-V_{t+1}^{PQ}(x+1)$, to approximate all $V_{t+1}(\mathbf{x})-V_{t+1}(\mathbf{x}+\mathbf{e}_i)$ for all $i\in\{1,2,...,m\}$, where $||\mathbf{x}||_1=x$. 
As a result, if the weight and volume distribution of different bookings vary significantly, the performance of PQ pricing may be negatively affected. Specifically, for bookings with extra small chargeable weight, such as type $i$ bookings, the term $(V_{t+1}^{PQ}(x)-V_{t+1}^{PQ}(x+1))/Q_i$ tends to be substantial, resulting in a significant increase in $r_i^{PQ}(t,x)$. Conversely, bookings with a very large chargeable weight will have significantly reduced pricing. In essence, extreme pricing issues may arise under PQ pricing strategy. 
To address this, we propose an augmented quantity-based approximation method. The main idea of this method is to replace $Q_i$ with $Q^s$ in the term $(V_{t+1}^{PQ}(x)-V_{t+1}^{PQ}(x+1))/Q_i$ to avoid extreme pricing, where $Q^s=\sum_{i=1}^m p_iQ_i$ denotes the weighted mean of $Q_i$ ($i=1,2,...,m$). 
Concretely, let $J_t^{AQ}(x)$ denote the expected revenue-to-go in our augmented method. At period $T$, $J_T^{AQ}(x)$ equals $V_T^{PQ}(x)$, as shown in Equation \eqref{eq:PQ_bound}:

\begin{equation*}
    J_T^{AQ}(x)=-E\left[h_w\left((\sum_{l=1}^xW_l^s-C_w)^+\right)+h_v\left((\sum_{l=1}^xV_l^s-C_v)^+\right)\right].
\label{eq:AQ_bound}
\end{equation*}

The augmented quantity-based (AQ) pricing for type $i$ bookings at state $x$ and period $t$, denoted as $r_i^{AQ}(t,x)$, is determined by solving the following equation:

\begin{equation*}
    r_i = \frac{1-F_t^i(r_i)}{f_t^i(r_i)}+\frac{J_{t+1}^{AQ}(x)-J_{t+1}^{AQ}(x+1)}{Q^s}.
\label{eq:AQ_price}
\end{equation*}

$J_t^{AQ}(x)$ is determined through backward induction using the following recursive equation:

\begin{equation*}
\begin{aligned}
    J_t^{AQ}(x)=&\sum_{i=1}^mb_t^i(r_i^{AQ}(t,x))\left(r_i^{AQ}(t,x)Q_i-(J_{t+1}^{AQ}(x)-J_{t+1}^{AQ}(x+1))\right) + J_{t+1}^{AQ}(x)  \\
    & \quad\quad\quad\quad\quad t=0,1,2,...,T-1.
\end{aligned}
\label{eq:JAQ}
\end{equation*}

In the general model framework, our augmented quantity-based approximation method prices type $i$ bookings at state $\mathbf{x}$ and period $t$ with $r_i^{AQ}(t,||\mathbf{x}||_1)$.

\section{Weight-Volume-based Approximation Method}
\label{sec:weight}
This section presents an approximation method based on weight and volume. 
By replacing the random weight and volume with their expected values, we create a new weight-volume-based MDP (WV-MDP). The state of the WV-MDP is $(w, v)$, where $w$ represents the cumulative expected weight of accepted bookings, and $v$ represents the cumulative expected volume of accepted bookings. We denote $\Tilde{V}_t(w, v)$ as the maximum expected revenue.
The Bellman equation and boundary condition for this MDP are as follows:

\begin{align}
    \Tilde{V}_t(w,v)=&\sum_{i=1}^m{\rm sup_{r_i}}\big\{b_t^i(r_i)(r_iQ_i-(\Tilde{V}_{t+1}(w,v)- \Tilde{V}_{t+1}(w+w_i,v+v_i)))\big\} + \Tilde{V}_{t+1}(w,v), \nonumber \\
    & \quad\quad\quad\quad\quad\quad t=0,1,2,...,T-1. \label{eq:wv} \\
    \Tilde{V}_T(w,v)=&-h_w(w-C_w)^+-h_v(v-C_v)^+ \label{eq:wv_boundary} 
\end{align}
where $w_i=EW_i,\ v_i=EV_i$.

\begin{figure}
  \centering
  \includegraphics[width=0.75\linewidth]{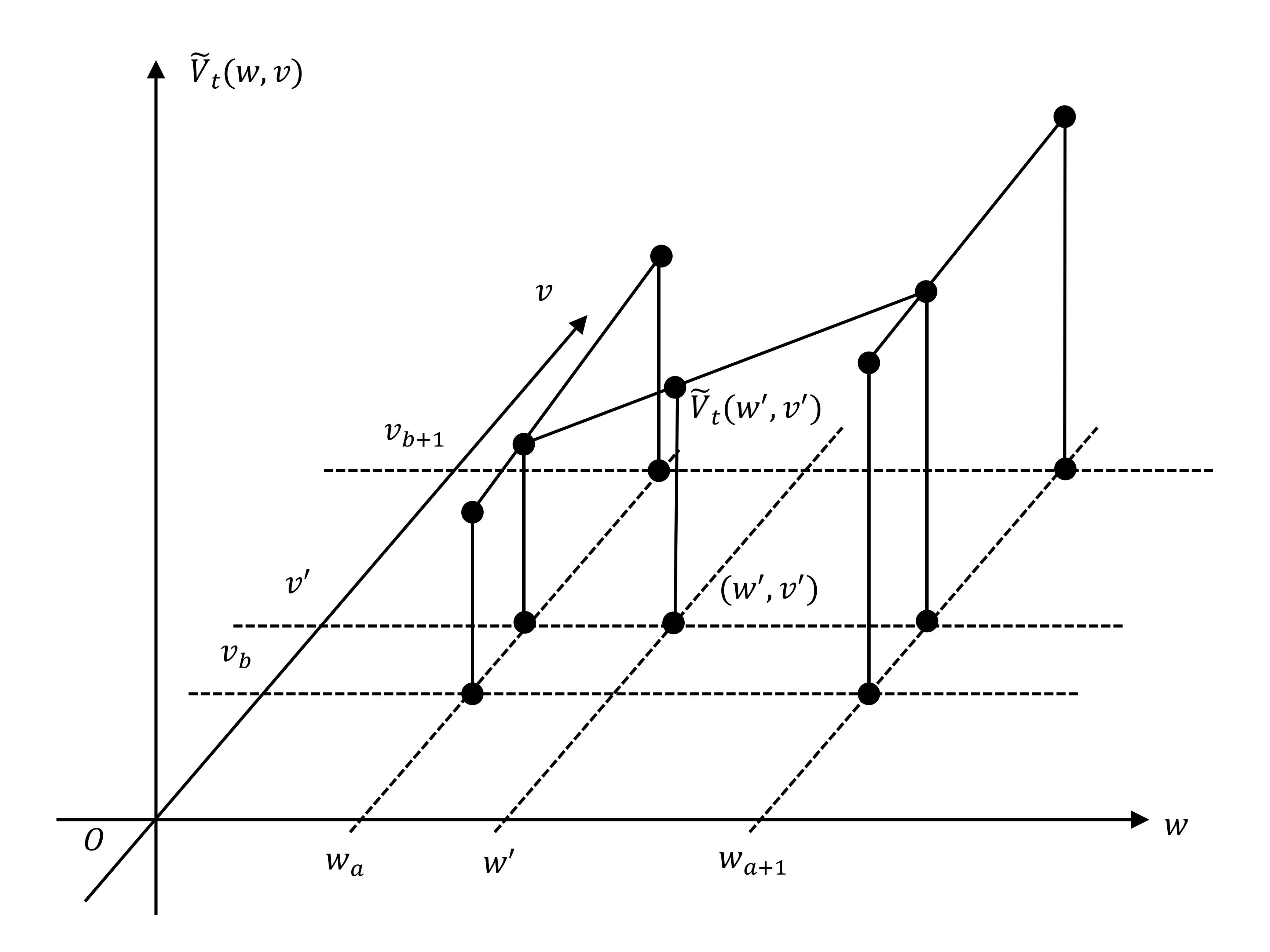}
  \caption{Bilinear interpolation}
  \label{fig:interpn}
\end{figure}

In the process of solving Equations \eqref{eq:wv}-\eqref{eq:wv_boundary} using backward induction, a discretization approach is applied to the continuous state variables $w$ and $v$. Specifically, the weight is divided into $A$ segments with an interval size of $\Delta w$, and the volume is divided into $B$ segments with an interval size of $\Delta v$. This division results in $K = (A+1)(B+1)$ sampling points, denoted as $(w_a, v_b)$ ($a = 0,1,\dots,A$ and $b = 0,1,\dots,B$). Starting from the boundary condition \eqref{eq:wv_boundary}, we determine the values of the sampling points $\Tilde{V}_t(w_a,v_b)$ based on Equation \eqref{eq:wv} and the function $\Tilde{V}_{t+1}(w,v)$. For non-sample points, the estimation of $\Tilde{V}_{t+1}(w',v')$ is accomplished through bilinear interpolation. Specifically, for the point $(w',v')$, we identify the smallest sampling grid that contains it, i.e., $w_a\leq w' \leq w_{a+1}$ and $v_b\leq v' \leq v_{b+1}$. The approximate value of $\Tilde{V}_{t+1}(w',v')$ is given by Equation \eqref{eq:interpn}, and it is illustrated in Fig. \ref{fig:interpn}.
\begin{equation}
\begin{aligned}
    \Tilde{V}_{t+1}(w',v') =& (1-\alpha)(1-\beta)\Tilde{V}_{t+1}(w_a,v_b) + (1-\alpha)\beta\Tilde{V}_{t+1}(w_a,v_{b+1})\\
    &+\alpha(1-\beta)\Tilde{V}_{t+1}(w_{a+1},v_b) + \alpha\beta\Tilde{V}_{t+1}(w_{a+1},v_{b+1}),
    \label{eq:interpn}
\end{aligned}
\end{equation}
where $\alpha=(w'-w_a)/\Delta w,\beta=(v'-v_b)/\Delta v$.

Next, we will discuss two pricing strategies based on $\Tilde{V}_t(w,v)$ in Sections \ref{subsec:PWVA} and \ref{subsec:WVAS}, respectively. The first strategy considers only the mean weight and volume of bookings, while the second strategy also incorporates second-order information.

\subsection{Primal Weight-Volume-based Approximation Method}
\label{subsec:PWVA}
Define 
\begin{equation}
V_t^{WV}(\mathbf{x})=\Tilde{V}_t(\sum_{i=1}^mx_iw_i,\sum_{i=1}^mx_iv_i).
\label{eq:def_wv}
\end{equation}
We employ the value of $V_t^{WV}(\mathbf{x})$ as an approximation for $V_t(\mathbf{x})$. 
Note that if we discretize the weight and volume with sufficient precision, there would be no need for interpolation to calculate $V_t^{WV}(\mathbf{x})$. This scenario is referred to as the ideal scenario.
The primal weight-volume-based (WV) approximation method utilizes the solution of the optimization problem:

\begin{equation}
    {\rm sup_{r_i}}\big\{m_{t}^i(1-F_t^i(r_i))\left(r_iQ_i-(V_t^{WV}(\mathbf{x})-V_t^{WV}(\mathbf{x}+\mathbf{e}_i))\right)\big\},
    \label{eq:opt_wv}
\end{equation}
as the pricing for a type $i$ booking at period $t$ and state $\mathbf{x}$. This pricing is denoted as $r_i^{WV}(t,\mathbf{x})$. Similar to Theorem \ref{prop:V_x}, WV pricing $r_i^{WV}(t,\mathbf{x})$ satisfies the following corollary:
\begin{corollary}
    A necessary condition for the price $r_i$ to be the solution of the optimization problem \eqref{eq:opt_wv} at period $t$ and state $\mathbf{x}$ is:

    \begin{equation}
        r_i = \frac{1-F_t^i(r_i)}{f_t^i(r_i)}+\frac{V_{t+1}^{WV}(\mathbf{x})-V_{t+1}^{WV}(\mathbf{x}+\mathbf{e}_i)}{Q_i}.
        \label{eq:ropt_wv}
    \end{equation}
    If the function $(1-F_t^i(r_i))^2/f_t^i(r_i)$ is decreasing in $r_i$, then the Equation \eqref{eq:ropt_wv} has a unique solution.
    \label{prop:Vwv_x} 
\end{corollary}

Next, we explore several structural properties of WV method. We first consider the Certainty Equivalent (CE) approach. The CE approach simplifies the treatment of uncertainty in a model via substituting random variables with their respective expectations \citep{bertsekas2012dynamic}.
Let $\bar{V}_t(\mathbf{x})$ represent the CE maximum expected revenue-to-go. Equations \eqref{eq:CE}-\eqref{eq:CE_boundary} below characterize the CE model.

\begin{align}
    \bar{V}_t(\mathbf{x})=&\sum_{i=1}^m{\rm sup_{r_i}}\big\{b_t^i(r_i)\left(r_iQ_i-(\Bar{V}_{t+1}(\mathbf{x})-\Bar{V}_{t+1}(\mathbf{x}+\mathbf{e_i}))\right)\big\} + \Bar{V}_{t+1}(\mathbf{x}),  \nonumber \\ 
    &\quad\quad\quad t=0,1,2,...,T-1. \label{eq:CE} \\
    \bar{V}_T(\mathbf{x})=&-h_w\left((\sum_{i=1}^{m}\sum_{j=1}^{x_i}E[W_{ij}]-C_w)^+\right)-h_v\left((\sum_{i=1}^{m}\sum_{j=1}^{x_i}E[V_{ij}]-C_v)^+\right). \label{eq:CE_boundary}
\end{align}

The CE model shares the same Bellman equation as the general model. The distinction lies in the boundary condition of the CE model, where the random variables in the boundary Equation \eqref{eq:gen_bound} of the general model are replaced by their expectations. Using backward induction, we can solve the CE model and compute the CE pricing, denoted as $\bar{r_i}(t,\mathbf{x})$. The following proposition demonstrates that $\bar{V}_t(\mathbf{x})$ effectively serves as an upper bound for $V_t(\mathbf{x})$.

\begin{proposition}
    $\bar{V}_t(\mathbf{x}) \geq V_t(\mathbf{x}), \quad \forall t=0,1,...,T.$
\label{prop:compare_1}
\end{proposition}

The next proposition states that $V_t^{WV}(\mathbf{x})$ equals $\bar{V}_t(\mathbf{x})$ in the ideal scenario.
\begin{proposition}
    In the ideal scenario, $V_t^{WV}(\mathbf{x})=\Bar{V}_t(\mathbf{x}),\ \forall \ t=0,1,...,T.$
\label{prop:idealwv}
\end{proposition}

Propositions \ref{prop:compare_1} and \ref{prop:idealwv} indicate that the limit of WV pricing is the CE pricing, and $V_t^{WV}(\mathbf{x})$ serves as an upper bound for $V_t(\mathbf{x})$ in the ideal scenario. Therefore, to analyze the gap in expected revenue between the WV pricing strategy in the ideal scenario and the optimal pricing strategy, we can compare the gap between the CE pricing strategy and the optimal pricing strategy. 
Let $J_t^{CE}(\mathbf{x})$ denote the expected revenue-to-go using the CE pricing strategy in the general model. According to the optimality principle, we can conclude that $V_t(\mathbf{x}) \geq J_t^{CE}(\mathbf{x})$. Combined with Proposition \ref{prop:compare_1}, it implies that  
\begin{equation}
    \Bar{V}_t(\mathbf{x}) \geq V_t(\mathbf{x}) \geq J_t^{CE}(\mathbf{x}).
    \label{eq:compare_2}
\end{equation}

Let $D_t^{CE}(\mathbf{x})=\Bar{V}_t(\mathbf{x})-J_t^{CE}(\mathbf{x})$. Equation \eqref{eq:compare_2} indicates $0\leq V_{t}(\mathbf{x})-J_t^{CE}(\mathbf{x}) \leq D_t^{CE}(\mathbf{x})$. Thus, the upper bound for $D_t^{CE}(\mathbf{x})$ also serves as an upper bound for $V_{t}(\mathbf{x})-J_t^{CE}(\mathbf{x})$. The following theorem presents an upper bound for $D_t^{CE}(\mathbf{x})$ when both functions $h_w$ and $h_v$ are linear.

\begin{theorem}
    When $h_w(x)=h_w x, h_v(x)=h_v x$, 
    \begin{equation*}
    \begin{aligned}
        D_t^{CE}(\mathbf{x})\leq \max\limits_{||\mathbf{x}||_1\leq x_{max}}&\Bigg\{
        \frac{h_w}{2}(\sqrt{\sigma_W^2(\mathbf{x})+(C_w-\Bar{w}(\mathbf{x}))^2}-|C_w-\Bar{w}(\mathbf{x})|)\\
        &+ \frac{h_v}{2}(\sqrt{\sigma_V^2(\mathbf{x})+(C_v-\Bar{v}(\mathbf{x}))^2}-|C_w-\Bar{v}(\mathbf{x})|)
        \Bigg\}, \\
        &\quad\quad\quad \forall ||\mathbf{x}||_1\leq x_{max},\ t=0,1,...,T,
    \end{aligned}
    \end{equation*}
    where $\Bar{w}(\mathbf{x}),\ \sigma_W^2(\mathbf{x})$ and $\Bar{v}(\mathbf{x}),\ \sigma_V^2(\mathbf{x})$ denote the mean and variance of 
    $\sum_{i=1}^{m}\sum_{j=1}^{x_i}W_{ij}$ and $\sum_{i=1}^{m}\sum_{j=1}^{x_i}V_{ij}$ respectively.
\label{prop:DCE}
\end{theorem}

\begin{remark}
    Define $M(\mathbf{x})\triangleq  0.5h_w\left(\sqrt{\sigma_W^2(\mathbf{x})+(C_w-\Bar{w}(\mathbf{x}))^2}-|C_w-\Bar{w}(\mathbf{x})|\right) 
        + 0.5h_v\left(\sqrt{\sigma_V^2(\mathbf{x})+(C_v-\Bar{v}(\mathbf{x}))^2}-|C_w-\Bar{v}(\mathbf{x})|\right).$
    Let $\mathbf{x}^* = \text{argmax}_{||\mathbf{x}||_1 \leq x_{\text{max}}} \{ M(\mathbf{x})\}$. By setting $t = 0$ and $\mathbf{x} = \mathbf{0}$ in Theorem \ref{prop:DCE}, the difference between the expected revenue of WV pricing strategy in the ideal scenario and the optimal expected revenue, denoted as $D_0^{\text{CE}}(\mathbf{0})$, satisfies:
     $D_0^{\text{CE}}(\mathbf{0}) \leq M(\mathbf{x}^*) \leq 0.5[h_w\sigma_W(\mathbf{x}^*) + h_v\sigma_V(\mathbf{x}^*)].$
     Intuitively, from this inequality, we can infer that the performance of WV pricing strategy deteriorates as the uncertainty of weight and volume of the bookings increases. Furthermore, this inference will be validated in the numerical experiments.
     When neither the weight nor the volume of the bookings is uncertain, the right side of the inequality equals 0, resulting in no difference in the expected revenue.
\end{remark}

\subsection{Weight-Volume-based Approximation Method with Second-Order Information}
\label{subsec:WVAS}

Considering the inherent uncertainty in the model, relying solely on the accumulated expected weight and volume to represent the state $\mathbf{x}$ at period $t$ in Equation \eqref{eq:def_wv} could introduce bias. In cases where uncertainty is significant, and the penalty for unloading cargo is severe, airlines tend to adopt a conservative and risk-averse approach. Consequently, the perceived accumulated weight and volume may exceed expectations. To address this, we incorporate second-order information and propose the weight-volume-based pricing strategy with second-order information (WVS pricing).
Concretely, define 

\begin{equation*}
    V_t^{WVS}(\mathbf{x})=\Tilde{V}_t(\sum_{i=1}^mx_i(w_i+\theta \sigma_{W_i}),\sum_{i=1}^mx_i(v_i+\theta \sigma_{V_i})),
\label{eq:def_WVS} 
\end{equation*}
where $\theta \in \mathbb{R}$ is a coefficient. We consider $w_i+\theta \sigma_{W_i}$ as the perceived weight and $v_i+\theta \sigma_{V_i}$ as the perceived volume of a type $i$ booking. Similar to WV method, we use the value of $V_t^{WVS}(\mathbf{x})$ as an approximation of $V_t(\mathbf{x})$. Subsequently, we can compute WVS pricing $r_i^{WVS}(\mathbf{x})$ in the same manner. 
The optimal value of $\theta$ can vary depending on various factors, such as the level of uncertainty, and the penalty associated with unloading cargo. Therefore, determining the optimal value is challenging. To address this issue, we employ a line search method in the numerical experiments to find a value of $\theta$ close to the optimal solution.

\section{Numerical Experiments}
\label{sec:NE}
Numerical experiments are undertaken to compare the performance of various proposed methods in this section. Two example sets are designed. The first set comprises three types of bookings, which allows the computation of optimal pricing and the comparison of different methods' performance against the optimal strategy. 
The second set includes twenty-seven types of bookings, with the settings based on \cite{Huang2010} and a real-world dataset collected by Megacap Aviation Service Corporation, a major aircargo forwarder in China. In this case, the exact solution of the MDP is not feasible. 
Hence, different approximation methods' performance are compared against an upper bound of the optimal expected revenue, $\bar{V}_0(\mathbf{0})$. 
It is assumed that the reservation prices for each type of booking at any given time follow a Weibull distribution, which allows for the modeling of diverse reservation price behaviors \citep{bitran1997periodic}.
The cumulative distribution function (CDF) of a Weibull distribution with scale parameter $\alpha$ and shape parameter $\beta$ is $F(x;\alpha,\beta)=1-e^{-(\frac{x}{\alpha})^\beta}$. 

We perform a three-factor test in each numerical set, with the factors including the ratio of capacity to demand, the penal factor on the boundary, and the coefficient of variation for weight and volume.
Concretely, let $D_w=\sum_{i=1}^m\int_0^L\lambda_s^iw_ids$ denote the total expected arrival weight, and $D_v=\sum_{i=1}^m\int_0^L\lambda_s^iv_ids$ denote the total expected arrival volume. The ratio of capacity to demand ($C/D$) is defined as $C/D=\frac{C_w}{D_w}=\frac{C_v}{D_v}$. We examine four distinct values, 0.8, 0.9, 1.0, and 1.1, for $C/D$. Let $\eta_w$ and $\eta_v$ denote the approximate revenue of a unit of weight and volume, respectively, and we utilize the following equation to calculate them. 

\begin{equation*}
    \sum_{i=1}^m\int_0^L\lambda_s^iE(RP_s^i)Q_ids=\eta_wD_w=\eta_vD_v.
\end{equation*}

We assume the penalty functions in the boundary condition are linear in excess capacity as follows.

\begin{equation*}
    h_w\left((x-C_w)^+\right) = pf\eta_w (x-C_w)^+, \label{eq:h_w}
\end{equation*}
\begin{equation*}
    h_v\left((x-C_v)^+\right) = pf\eta_v (x-C_v)^+, \label{eq:h_v} 
\end{equation*}
where $pf\in \mathbb{R}$ is the penal factor, for which we consider three different values: 1, 1.25, and 1.5. 
All weight and volume distributions for each booking type are assumed to follow independent normal distributions with the same coefficient of variation ($cv$), i.e., $cv=\sigma_{W_i}/w_i=\sigma_{V_i}/v_i, \ \forall i\in\{1,2,...,m\}$. Also, we examine three distinct values, 0.2, 0.3, and 0.5, for $cv$ to represent varying uncertain degrees. Hence, in each numerical experiment, we assess $4\times 3 \times 3 =36$ scenarios. 

\subsection{Numerical Results} \label{sucsec:NumRes}
\subsubsection{A toy Example of three types of bookings.}
\paragraph{}
We first consider a set of three types of bookings, i.e., $m=3$. Table \ref{tab:wv_3} displays the mean weight and volume for each type of booking.
\begin{table}
  \centering
  \caption{Mean weight and volume of each type of booking}
  \begin{tabular}{lrrr}
    \hline
    type &1&2&3 \\
    \hline
    $w_i(\text{kg})$ & 100 & 75 & 150  \\
    $v_i(\times 10^4 \text{cm}^3)$ & 60 & 50 & 75   \\
    \hline
  \end{tabular}
  \label{tab:wv_3}
\end{table}
The expected densities of these three types of bookings are equal to, less than, and more than the inverse of $\gamma$, respectively.

We set $L=75$. The arrival process for type $i$ bookings corresponds to a non-homogeneous Poisson process $PP(\lambda^i_s)$. The reservation price distribution for type $i$ bookings at time $s$ follows Weibull$(\alpha^i_s, \beta^i_s)$, $s\in[0,L]$. Table \ref{tab:para_3_nonhomo} shows the specific parameter values.
\begin{table}
  \centering
  \caption{Parameters for each type of booking}
  \begin{tabular}{llcc}
    \hline
    type&  \multicolumn{1}{c}{$ \lambda^i_s$} & $\alpha^i_s$ & $\beta^i_s$ \\
    \hline
    1 & $\left\{
    \begin{array}{ll}
    0.12s/L+0.04, &  \quad 0\leq s\leq 2L/3; \\
    -0.12s/L+0.2,  &  \quad 2L/3\leq s \leq L. \\
    \end{array} \right.$ & $4+2s/L$  & 5  \\
    2 & $\left\{
    \begin{array}{ll}
    0.075s/L+0.025, &  \quad 0\leq s\leq 2L/3; \\
    -0.075s/L+0.125,  & \quad 2L/3\leq s \leq L. \\
    \end{array} \right.$ & $3+1.5s/L$ & 5  \\
    3  & $\left\{
    \begin{array}{ll}
    0.09s/L+0.03, &  \quad 0\leq s\leq 2L/3; \\
    -0.09s/L+0.15,  & \quad 2L/3 \leq s \leq L. \\
    \end{array} \right.$ & $3+1.5s/L$ & 5   \\
    \hline
  \end{tabular}
  \label{tab:para_3_nonhomo}
\end{table}
In this setting, the arrival rate for each type of booking initially increases linearly and then decreases linearly. The total expected arrival weight, $D_w$, is 1690 kg, and the volume, $D_v$, is 960 $\times 10^4 \ \text{cm}^3$.
The expected reservation price for each type of booking increases linearly, eventually reaching 1.5 times its initial value. This pattern aligns with the common understanding that flight prices tend to increase as the departure time approaches.
Let $\Delta t = 1/3$ and $T = L/\Delta t = 225$. Then the probability of at most one booking arriving during a period is less than 99\%. 

In this case, the general model can be solved to get the optimal expected revenue. Denote the percentage difference of the expected revenue for an approximation method from the optimal expected revenue as follows:

\begin{equation*}
    \Delta^h V = \frac{V_0(\mathbf{0})-J_0^{h}(\mathbf{0})}{V_0(\mathbf{0})}\times 100\%,
    \label{eq:diff}
\end{equation*}
where $J_0^{h}(\mathbf{0})$ denotes the expected revenue for pricing strategy $h$, $h\in \{\text{PQ,AQ,WV,WVS}\}$. 

When solving WV(WVS) pricing, we divide the weight into $50$ segments, each with a length of
$50$ kg, and the volume into $50$ segments, each with a length of $30\times10^4 \ \text{cm}^3$. The results of $\Delta^h V$ and the optimal $\theta^*$ (at a granularity of 0.01) for WVS pricing are recorded in Table \ref{tab:res_3}. The average computation time for the quantity-based methods is approximately 5 seconds, while for the weight-volume-based methods with a specific $\theta$, it is around 2 minutes. 

\begin{table}
    \centering
    \caption{Percentage difference from the optimal expected revenue for three types of bookings}
    \resizebox{\textwidth}{!}{
    \begin{tabular}{llcccccccccccc} 
    \hline
        $pf$ & $C/D$ & \multicolumn{4}{c}{$cv=0.2$}  & \multicolumn{4}{c}{$cv=0.3$} & \multicolumn{4}{c}{$cv=0.5$} \\ 
        \cmidrule(lr){3-6} \cmidrule(lr){7-10} \cmidrule(lr){11-14}   
        ~ & ~ & PQ & AQ & WV & WVS ($\theta^*$) & PQ & AQ & WV & WVS ($\theta^*$) & PQ & AQ & WV & WVS ($\theta^*$) \\ 
        \hline
        1 & 0.8 & 1.49 & 0.66 & 0.11 & 0.09 (0.04) & 1.41 & 0.53 & 0.29 & 0.22 (0.05) & 1.42 & 0.40 & 0.88 & 0.58 (0.09) \\ 
        ~ & 0.9 & 0.94 & 0.48 & 0.07 & 0.06 (0.04) & 0.89 & 0.39 & 0.21 & 0.16 (0.05) & 0.89 & 0.30 & 0.70 & 0.44 (0.09) \\ 
        ~ & 1 & 0.56 & 0.32 & 0.06 & 0.05 (0.04) & 0.53 & 0.27 & 0.16 & 0.12 (0.06) & 0.53 & 0.21 & 0.54 & 0.34 (0.09) \\ 
        ~ & 1.1 & 0.31 & 0.20 & 0.04 & 0.03 (0.04) & 0.29 & 0.16 & 0.11 & 0.08 (0.06) & 0.30 & 0.14 & 0.37 & 0.23 (0.09) \\ 
        1.25 & 0.8 & 1.75 & 0.83 & 0.16 & 0.09 (0.06) & 1.67 & 0.68 & 0.44 & 0.23 (0.09) & 1.71 & 0.52 & 1.42 & 0.61 (0.13) \\ 
        ~ & 0.9 & 1.14 & 0.61 & 0.11 & 0.07 (0.06) & 1.07 & 0.50 & 0.33 & 0.17 (0.09) & 1.09 & 0.39 & 1.12 & 0.48 (0.12) \\ 
        ~ & 1 & 0.69 & 0.41 & 0.09 & 0.05 (0.06) & 0.65 & 0.34 & 0.26 & 0.13 (0.09) & 0.67 & 0.28 & 0.86 & 0.38 (0.13) \\ 
        ~ & 1.1 & 0.39 & 0.25 & 0.06 & 0.04 (0.06) & 0.37 & 0.22 & 0.17 & 0.09 (0.09) & 0.39 & 0.18 & 0.58 & 0.27 (0.13) \\ 
        1.5 & 0.8 & 1.94 & 0.97 & 0.23 & 0.10 (0.08) & 1.89 & 0.80 & 0.65 & 0.24 (0.11) & 1.96 & 0.62 & 2.09 & 0.63 (0.15) \\ 
        ~ & 0.9 & 1.33 & 0.73 & 0.16 & 0.07 (0.08) & 1.23 & 0.60 & 0.49 & 0.18 (0.12) & 1.27 & 0.47 & 1.63 & 0.51 (0.16) \\ 
        ~ & 1 & 0.78 & 0.49 & 0.13 & 0.05 (0.09)  & 0.75 & 0.41 & 0.37 & 0.14 (0.12)  & 0.79 & 0.34 & 1.24 & 0.41 (0.16) \\ 
        ~ & 1.1 & 0.46 & 0.31 & 0.09 & 0.04 (0.09) & 0.43 & 0.26 & 0.24 & 0.09 (0.12) & 0.46 & 0.22 & 0.83 & 0.29 (0.15) \\ 
        ~ & ~ & ~ & ~ & ~  & ~ & ~ & ~ & ~ & ~ & ~ & ~ & ~ & ~  \\ 
        \multicolumn{2}{l}{Minimum} & 0.31 & 0.20 & 0.04 & 0.03  & 0.29 & 0.16 & 0.11 & 0.08 & 0.30 & 0.14 & 0.37 & 0.23 \\ 
        \multicolumn{2}{l}{Mean} & 0.98 & 0.52 & 0.11 & 0.06   & 0.93 & 0.43 & 0.31 & 0.15 & 0.96 & 0.34 & 1.02 & 0.43  \\ 
        \multicolumn{2}{l}{Maximum}  & 1.94 & 0.97 & 0.23 & 0.10 & 1.89 & 0.80 & 0.65 & 0.24 & 1.96 & 0.62 & 2.09 & 0.63 \\ \hline
    \end{tabular}
    }
    \label{tab:res_3}
\end{table}

In all 36 scenarios, the maximum percentage gaps of PQ, AQ, WV, and WVS pricing methods are all less than 2.1\%. The average percentage gap of PQ pricing is less than 1\%. In comparison, as an augmented variant of PQ pricing, AQ pricing consistently performs better, demonstrating an average improvement of nearly 0.5\% in terms of performance gap. When $cv$ is small, such as 0.2, WV pricing performs well, with a mean performance gap of 0.11\%, which aligns with Theorem \ref{prop:DCE}. WVS pricing incorporates second-order information into WV pricing to enhance its performance. The improvement from WVS pricing over WV pricing increases with higher $cv$ values. Specifically, when $cv=0.5$, the average improvement is 0.6\%, compared to 0.05\% when $cv=0.2$. In contrast, the quantity-based methods maintain stable performance as $cv$ increases. Interestingly, the performance of AQ pricing improves as $cv$ increases.

In general, as $pf$ increases, the performances of PQ, AQ, WV, and WVS pricing methods all deteriorate. However, the impact on the WVS pricing is minimal because the value of $\theta^*$ is adjusted to accommodate the increased impact of $pf$. Specifically, we observe that $\theta^*$ increases as both $cv$ and $pf$ increase. This observation aligns with intuition, as higher values of $cv$ and $pf$ lead to a more conservative approach in sales, resulting in an increase in the perceived weight and volume. Furthermore, the performance gap of all methods increases as the weight and volume capacity tightens ($C/D$ decreases). This observation is also intuitive because when weight and volume capacity are sufficient, the pricing methods become less valuable.

Overall, in the case of $m=3$ where there is low heterogeneity in the weight and volume distributions of diverse booking types, the quantity-based methods require less computation and yield great performance. Moreover, AQ pricing generates higher expected revenue than PQ pricing, while PQ pricing has more desirable structural properties as discussed in Section \ref{subsec:PQA}. The weight-volume-based methods are computationally more intensive, but demonstrate superior performance in scenarios characterized by low uncertainty in the weight and volume of bookings. The performance improvement from WVS pricing over WV pricing is relatively significant when both penalty costs and uncertainty levels are high. In addition, the quantity-based methods exhibit more robust performance than the weight-volume-based methods in high uncertainty scenarios.

\subsubsection{An Example of 27 types of bookings.}
\paragraph{}
Now we consider a set of 27 types of bookings, i.e., $m=27$. The dataset comprises cargo booking information for flights from PVG Airport to DEL Airport between January $2^{nd}$, 2019 and May $24^{th}$, 2021. It encompasses over 10,000 records, including the weight and volume of bookings, the booking placement time, flight numbers, and flight departure time.

Suppose that there are nine cargo categories with different mean weights and volumes shown in Table \ref{tab:wv_27} \citep{Huang2010}.
\begin{table}[b]
  \centering
  \caption{Mean weight and volume for cargo categories}
  \begin{tabular}{lrrrrrrrrr}
    \hline
    category &1&2&3&4&5&6&7&8&9 \\
    \hline
    $w_i(\times 10 \text{kg})$ & 8 & 16 & 40 & 10 & 20 & 50 & 10 & 20 & 50  \\
    $v_i(\times 10^5 \text{cm}^3)$ & 6 & 12 & 30 & 6 & 12 & 30 & 5 & 10 & 25   \\
    \hline
  \end{tabular}
  \label{tab:wv_27}
\end{table}
Categories 1-3 are loose cargo, given that their volume-to-weight ratios surpass standard $\gamma$ ($=0.6 \times 10^4\text{cm}^3/\text{kg}$). Conversely, categories 7-9 are dense cargo due to their lower volume-to-weight ratio. Categories 4-6 are normal cargo with their volume-to-weight ratio equal to $\gamma$. 
Table \ref{tab:p_cat} shows the probability distribution for the categories of bookings.
\begin{table}[b]
  \centering
  \caption{Probability distribution for the categories}
    \begin{tabular}{lccccccccc}
    \hline Category & 1 & 2 & 3 & 4 & 5 & 6 & 7 & 8 & 9 \\
    \hline Probability & 0.0833 & 0.0833 & 0.0833 & 0.1667 & 0.1668 & 0.1667 & 0.0833 & 0.0833 & 0.0833 \\
    \hline
    \end{tabular}
  \label{tab:p_cat}
\end{table}

\begin{figure}
  \centering
  \includegraphics[width=0.75\linewidth]{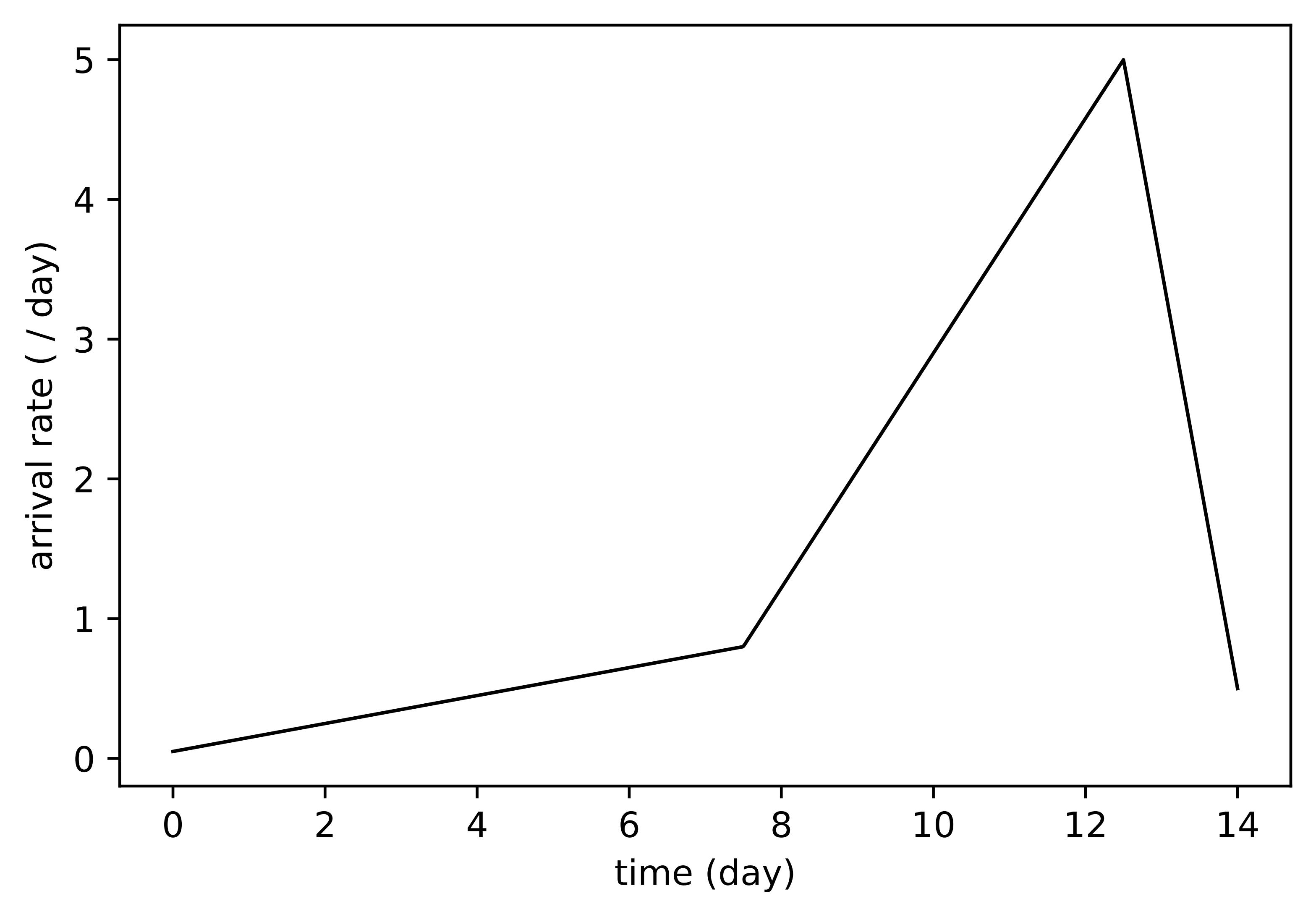}
  \caption{Total arrival rate}
  \label{fig:rate}
\end{figure}

We fit the arrival process using the real-world dataset and set $L=14$ days. The arrival process for all bookings follows a non-homogeneous Poisson process $PP(\lambda_s), \ s\in[0,L]$, where 

\begin{equation*}
    \lambda_s=\left\{
    \begin{array}{ll}
    0.1s+0.05, &  0\leq s\leq 7.5; \\
    0.8+0.84(s-7.5),  &  7.5\leq s \leq 12.5; \\
    2-3(s-13.5), &  12.5 \leq s \leq 14. \\
    \end{array} \right.
  \label{eq:rate_12}
\end{equation*}
The graph of arrival rate is depicted in Figure \ref{fig:rate}. It shows that the arrival rate initially increases slowly during the spot sale and rapidly increases midway through the period before diminishing as the flight approaches departure. In this setting, the expected total number of arrival bookings is $\int_0^L\lambda_sds \approx 21.8$. The total expected arrival weight $D_w$ and volume $D_v$ equal to 5530 kg and 3340 $\times 10^4 \ \text{cm}^3$, respectively. 
Let $\Delta t = 0.04$ and $T = L/\Delta t = 350$. Then the probability of at most one booking arriving during a period is near 99\%. 

In this numerical experiment, we categorize the reservation price distributions for a weight-volume distribution into three classes: median, high, and low reservation price distributions. Specifically, the parameters of the median reservation price distributions (Weibull($\alpha_s,\beta_s$)) for the nine categories of bookings are listed in Table \ref{tab:Weibull_M}.
\begin{table}
  \centering
  \caption{Parameters of the median reservation price distributions}
    \begin{tabular}{lccc}
    \hline Category & $1,4,7$ & $2,5,8$ & $3,6,9$  \\
    \hline $\alpha_s$ & $4+2s/L$ & $3.6+1.8s/L$ & $3.2+1.6s/L$ \\
            $\beta_s$ &  5  &  5  &  5 \\
    \hline
    \end{tabular}
  \label{tab:Weibull_M}
\end{table}
Additionally, for the high reservation price distribution, we set the scale parameter $\alpha_s$ to be 1.25 times that of the median. At the same time, for the low reservation price distribution, we set it to be 0.75 times that of the median. The probabilities for the medium, high, and low reservation price distributions are set at 0.4, 0.3, and 0.3 at any given time. Combining categories and reservation price distribution classes results in 27 types of bookings. For instance, the arrival rate of a booking with category 1 and a high reservation price distribution is $0.0833\times0.4\lambda_s$, and its reservation price distribution at time $s$ follows Weibull($1.25\times(4+2s/L)$,\ 5), $s\in[0,L]$.
Similar to the toy example, the expected reservation price for each type of booking increases linearly, eventually reaching 1.5 times its initial value.

In this case, because of the curse of dimensionality, solving the general model exactly is not feasible. Therefore, the expected revenues of different pricing methods are estimated using simulation. The 27 types of bookings' arrival is generated randomly based on the settings above. We run 5,000 simulations for each scenario and compute the average revenue from these simulations as the expected revenue of each pricing method. Simulation error is less than 0.5\% in almost all scenarios. By Equation \eqref{eq:def_wv} and Propositions \ref{prop:compare_1}-\ref{prop:idealwv}, we have 

\begin{equation*}
    V_t(\mathbf{0})\leq \bar{V}_t(\mathbf{0})=V_t^{WV}(\mathbf{0})=\Tilde{V}_t(0,0).
\end{equation*}

Since the exact value of $V_t(\mathbf{0})$ cannot be computed in this case, we calculate the relative percentage gaps of the expected revenue for an approximation method from its upper bound, $\Tilde{V}_t(0,0)$. The relative percentage differences are defined as follows:

\begin{equation*}
    \Delta^h \bar{V} = \frac{\Tilde{V}_0(0,0)-J_0^{h}(\mathbf{0})}{\Tilde{V}_0(0,0)}\times 100\%,
    \label{eq:diff_up}
\end{equation*}
where $h\in\{\text{PQ,AQ,WV,WVS}\}$.

When solving WV(WVS) pricing, we divide the weight into 50 segments, each with a length of $160 \ \text{kg}$, and the volume into 50 segments, each with a length of $100\times 10^4 \ \text{cm}^3$. The results of $\Delta^h \bar{V}$ and the reference optimal $\theta^*$ are recorded in Table \ref{tab:res_27}. Since the expected revenue for WVS pricing remains relatively stable near the optimal $\theta^*$ and considering the potential simulation error, the optimal $\theta^*$ recorded in Table \ref{tab:res_27} is a reference value, and there may be negative relative differences. In Section \ref{sec:insights}, we will further explore the performance improvements of WVS pricing compared to WV pricing and the correlations between the optimal value of $\theta$ and factors $C/D$, $pf$, and $cv$. The average computation time for the quantity-based methods is less than 1 minute, while for the weight-volume-based methods with one specific $\theta$, it amounts to around 2 hours.  

\begin{table}
    \centering
    \caption{Percentage difference from an upper bound for 27 types of bookings}
    \resizebox{\textwidth}{!}{
    \begin{tabular}{llrrrrrrrrrrrr}
    \hline
        $pf$ & $C/D$ & \multicolumn{4}{c}{$cv=0.2$}  & \multicolumn{4}{c}{$cv=0.3$} & \multicolumn{4}{c}{$cv=0.5$} \\ 
        \cmidrule(lr){3-6} \cmidrule(lr){7-10} \cmidrule(lr){11-14}   
        ~ & ~ & PQ & AQ & WV & WVS ($\theta^*$) & PQ & AQ & WV & WVS ($\theta^*$) & PQ & AQ & WV & WVS ($\theta^*$) \\ 
        1 & 0.8 & 8.89 & 3.53 & 0.69 & 0.68 (0.05) & 9.82 & 3.87 & 1.89 & 1.84 (0.05) & 12.85 & 6.21 & 5.15 & 4.17 (0.14) \\ 
        ~ & 0.9 & 5.07 & 2.23 & 0.58 & 0.54 (0.04) & 5.68 & 2.46 & 1.45 & 1.41 (0.06) & 9.04 & 4.63 & 4.14 & 3.45 (0.14)  \\ 
        ~ & 1 & 3.44 & 1.79 & 0.37 & 0.29 (0.07) & 4.11 & 2.51 & 0.67 & 0.56 (0.08) & 5.16 & 2.98 & 2.56 & 2.04 (0.17)  \\ 
        ~ & 1.1 & 1.59 & 1.00 & -0.16 & -0.18 (0.09) & 2.15 & 1.36 & 0.46 & 0.43 (0.06) & 3.37 & 2.23 & 1.83 & 1.38 (0.18) \\ 
        1.25 & 0.8 & 11.39 & 4.43 & 0.40 & 0.30 (0.04) & 12.63 & 5.52 & 1.88 & 1.64 (0.09) & 15.99 & 7.56 & 6.68 & 5.01 (0.19)  \\ 
        ~ & 0.9 & 6.42 & 3.06 & 1.01 & 0.91 (0.06) & 7.63 & 4.01 & 2.15 & 2.00 (0.09) & 10.15 & 5.48 & 5.23 & 4.03 (0.17)  \\ 
        ~ & 1 & 3.79 & 2.03 & -0.05 & -0.11 (0.07) & 4.76 & 3.11 & 1.41 & 1.15 (0.11) & 6.85 & 3.93 & 3.46 & 2.39 (0.18) \\ 
        ~ & 1.1 & 2.28 & 1.66 & 0.66 & 0.63 (0.04) & 1.71 & 0.95 & 0.42 & 0.31 (0.09) & 3.62 & 2.16 & 1.85 & 1.12 (0.16) \\ 
        1.5 & 0.8 & 12.95 & 5.59 & 0.73 & 0.57 (0.07) & 14.40 & 6.80 & 2.76 & 2.33 (0.09)  & 19.09 & 9.01 & 8.02 & 5.67 (0.23)  \\ 
        ~ & 0.9 & 8.20 & 3.69 & 1.20 & 1.11 (0.08) & 9.23 & 4.72 & 1.90 & 1.52 (0.12) & 12.22 & 6.36 & 6.78 & 4.80 (0.22)  \\ 
        ~ & 1 & 4.74 & 2.80 & 0.70 & 0.59 (0.08) & 5.77 & 3.78 & 1.52 & 1.30 (0.13) & 7.97 & 4.69 & 4.50 & 3.05 (0.21)  \\ 
        ~ & 1.1 & 2.29 & 1.48 & 0.54 & 0.47 (0.10) & 2.80 & 1.85 & 0.63 & 0.53 (0.12) & 4.51 & 2.99 & 3.27 & 2.21 (0.19)  \\ 
        ~ & ~ & ~ & ~ & ~ & ~ & ~ & ~ & ~ & ~ & ~ & ~ & ~ &    \\ 
        \multicolumn{2}{l}{Minimum} & 1.59 & 1.00 & -0.16 & \multicolumn{1}{c}{-0.18}  & 1.71 & 0.95 & 0.42 & \multicolumn{1}{c}{0.31}  & 3.37 & 2.16 & 1.83 & \multicolumn{1}{c}{1.12}  \\ 
        \multicolumn{2}{l}{Mean} & 5.92 & 2.77 & 0.56 & \multicolumn{1}{c}{0.48} & 6.72 & 3.41 & 1.43 & \multicolumn{1}{c}{1.25}  & 9.23 & 4.85 & 4.46 & \multicolumn{1}{c}{3.28}   \\ 
        \multicolumn{2}{l}{Maximum} & 12.95 & 5.59 & 1.20 & \multicolumn{1}{c}{1.11}  & 14.40 & 6.80 & 2.76 & \multicolumn{1}{c}{2.33} & 19.09 & 9.01 & 8.02 & \multicolumn{1}{c}{5.67}  \\ \hline
    \end{tabular}
    }
    \label{tab:res_27}
\end{table}

In all 36 scenarios, the average percentage gaps of PQ, AQ, WV, and WVS pricing methods are all less than 10\%. WV pricing demonstrates great performance in all test scenarios, particularly in scenarios with low uncertainty in the weight and volume of bookings. For instance, when $cv=0.2$, the mean performance gap of WV pricing method from the upper bound is only 0.56\%. WVS pricing method, which incorporates second-order information, further enhances the performance of the WV pricing method and performs best in this test set. Notably, the weight-volume-based pricing methods outperform the quantity-based methods, albeit requiring more computational efforts.

In general, the performances of PQ, AQ, WV, and WVS pricing methods all deteriorate as $C/D$ decreases, $pf$ increases, and $cv$ increases. However, the impact of $pf$ and $cv$ on the performances of AQ and WVS pricing are relatively lower compared to other methods. Our results show that the relative performances of these methods are similar to the case of $m=3$ above. However, as the heterogeneity of bookings increases, the performance gaps from the optimal expected revenue for all pricing methods are somewhat affected.

\subsection{Insights into Second-Order Information}\label{sec:insights}
Now we explore the performance improvements using WVS pricing compared to WV pricing and the correlations between the optimal value of $\theta$ and factors $C/D$, $pf$, and $cv$. 

\begin{figure}[b]
  \centering
  \subcaptionbox{$m=3$\label{fig:cd-m3}}
    {\includegraphics[width=0.49\linewidth]{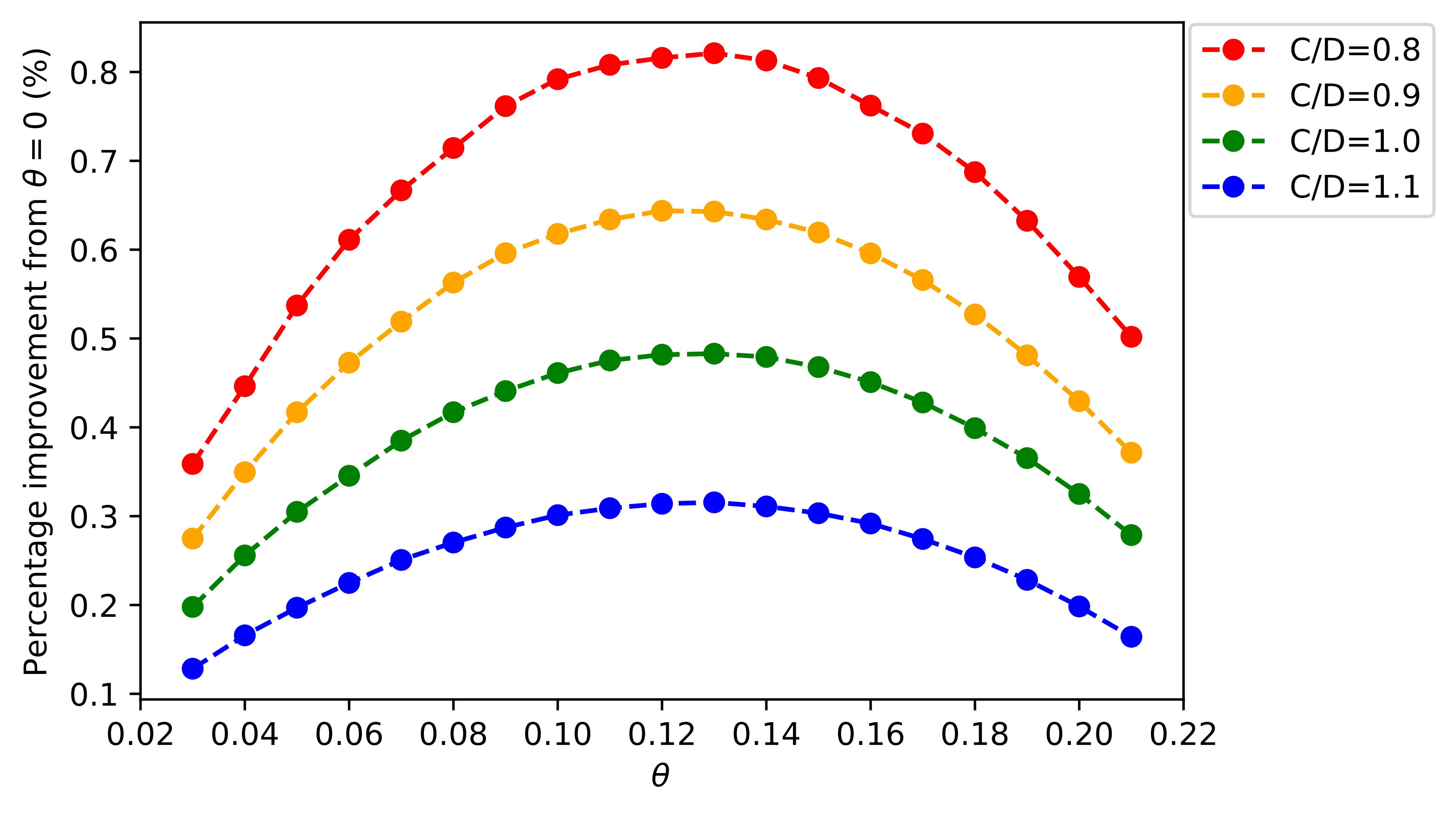}}
  \subcaptionbox{$m=27$\label{fig:cd-m27}}
    {\includegraphics[width=0.49\linewidth]{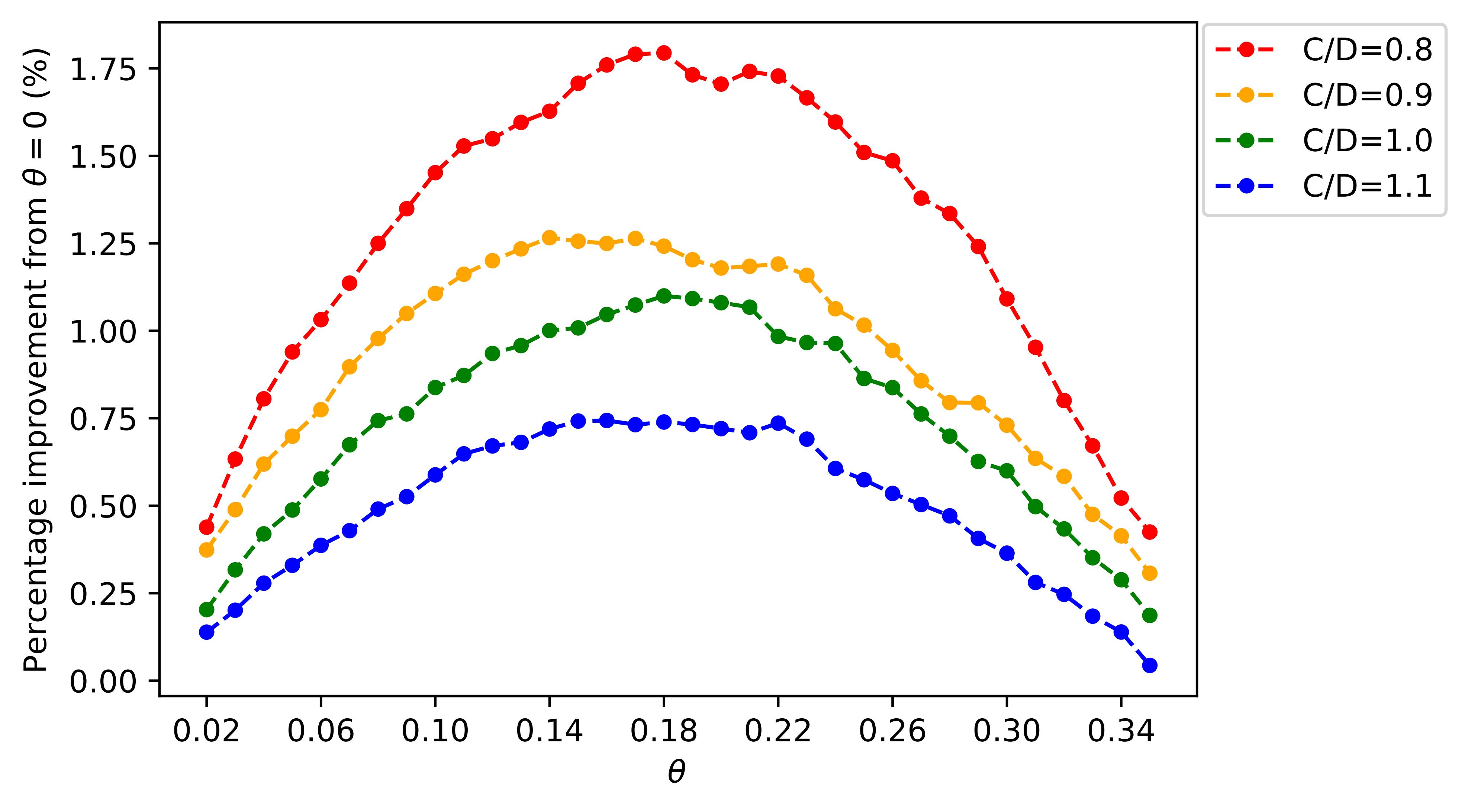}}
  \caption{Percentage improvements from WV pricing v.s. $\theta$ ($pf=1.25$, $cv=0.5$)}
  \label{fig:cd}
\end{figure}

We first discuss the impact of $C/D$. We vary $C/D\in\{0.8,0.9,1.0,1.1\}$, while keeping $pf=1.25$ and $cv=0.5$ fixed, following the tests in Section \ref{sucsec:NumRes}. Figures \ref{fig:cd-m3} and \ref{fig:cd-m27} show the performance improvements over WV pricing ($\theta=0$) with respect to $\theta$ in the cases of $m=3$ and $m=27$, respectively. (Note that the curves are not smooth in the case of $m=27$ due to simulation errors). 
The figures show that the percentage improvements increase as the capacity tightens, i.e., $C/D$ decreases. Moreover, it is observed that the improvements remain relatively stable near the optimal value of $\theta$. The optimal value of $\theta$ exhibits little correlation with $C/D$. In addition, the level of performance improvements in the case of $m=27$ is higher than that of $m=3$.

\begin{figure}
  \centering
  \subcaptionbox{$m=3$\label{fig:pf-m3}}
    {\includegraphics[width=0.49\linewidth]{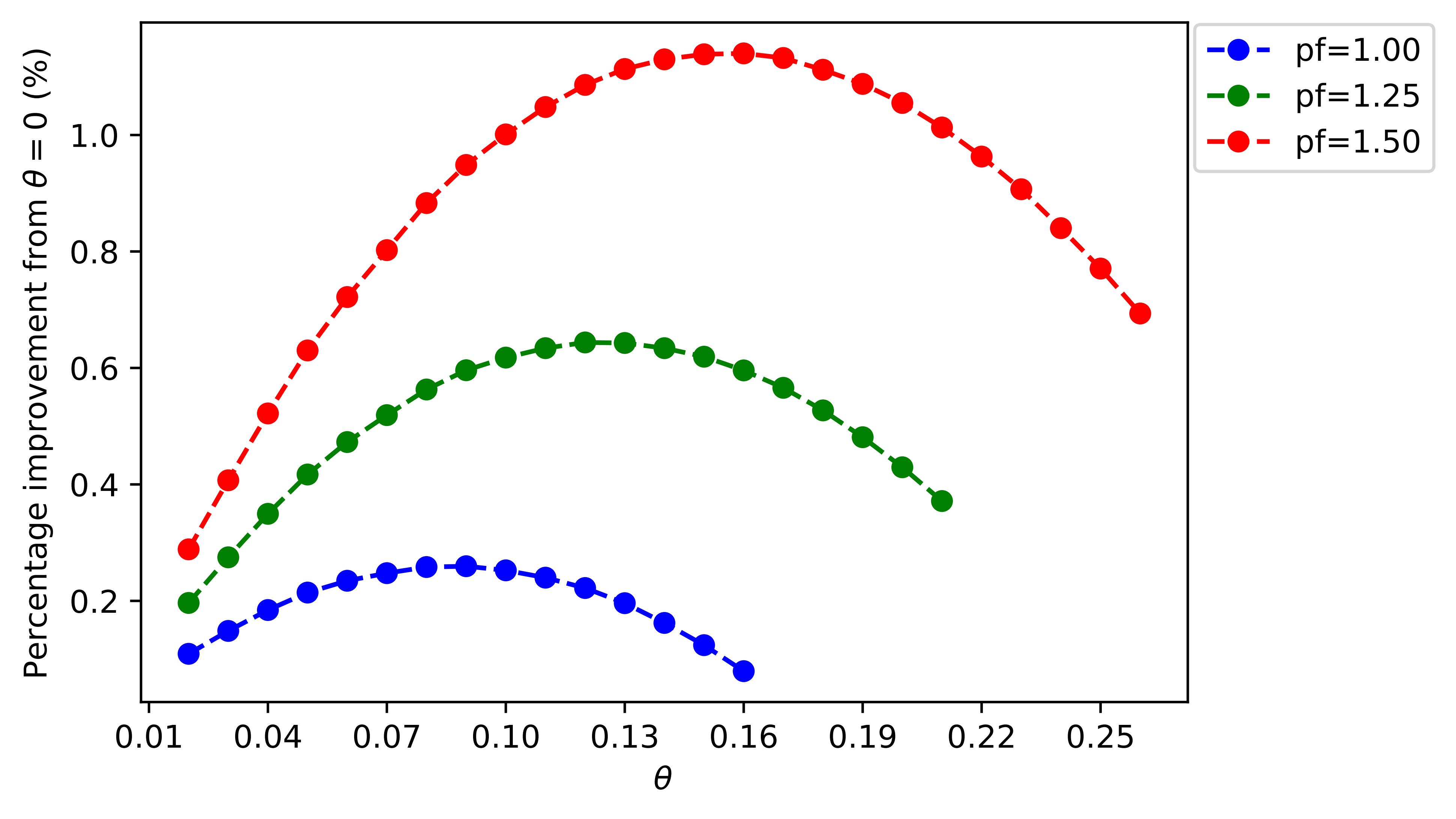}}
  \subcaptionbox{$m=27$\label{fig:pf-m27}}
    {\includegraphics[width=0.49\linewidth]{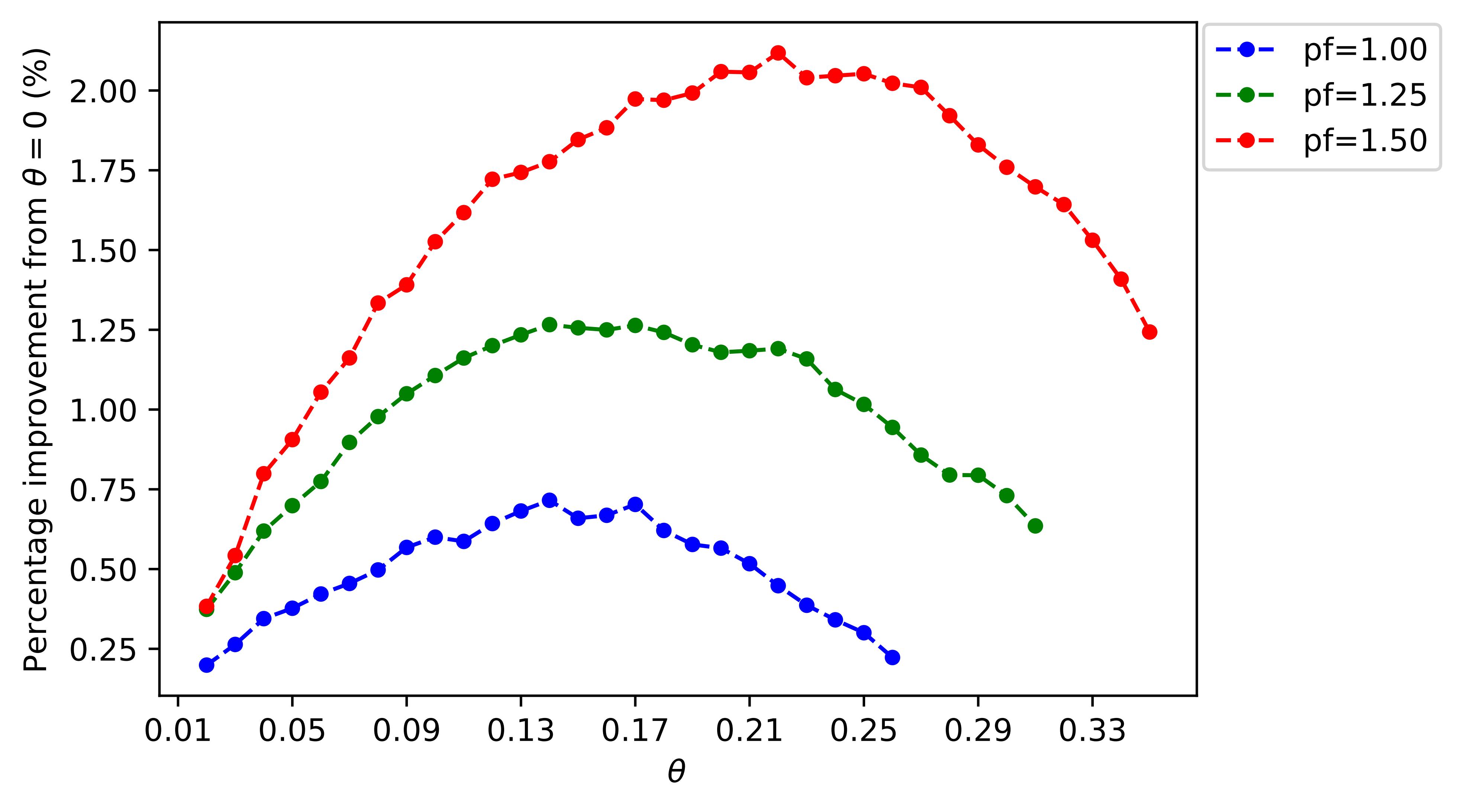}}
  \caption{Percentage improvements from WV pricing v.s. $\theta$ ($C/D=0.9$, $cv=0.5$)}
  \label{fig:pf}
\end{figure}

Next, we discuss the impact of $pf$. We vary $pf \in \{1.00, 1.25, 1.50\}$ while keeping $C/D=0.9$ and $cv=0.5$ fixed, following the tests in Section \ref{sucsec:NumRes}, as well. Figures \ref{fig:pf-m3} and \ref{fig:pf-m27} shows the performance improvements over WV pricing ($\theta=0$) with respect to $\theta$ in the cases of $m=3$ and $m=27$, respectively. 
The figures show that the percentage improvements increase as $pf$ increases. Moreover, the optimal value of $\theta$ also increases with higher $pf$ values. Consistent with the results in the discussion of the impact of $C/D$, the improvements remain relatively stable near the optimal value of $\theta$, and the level of performance improvements in the case of $m=27$ is more pronounced than that of $m=3$.

\begin{figure}
  \centering
  \subcaptionbox{$m=3$\label{fig:cv-m3}}
    {\includegraphics[width=0.49\linewidth]{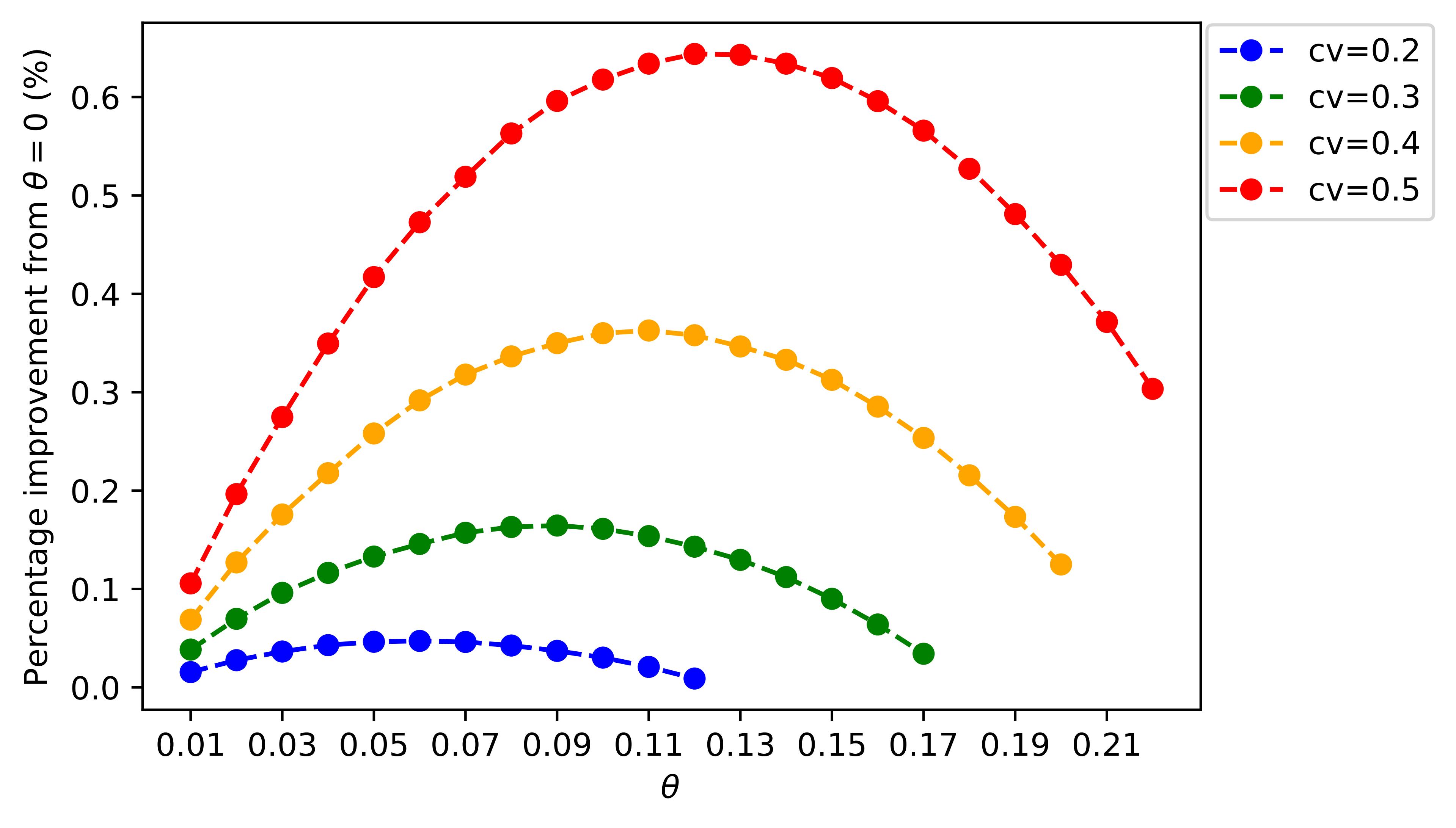}}
  \subcaptionbox{$m=27$\label{fig:cv-m27}}
    {\includegraphics[width=0.49\linewidth]{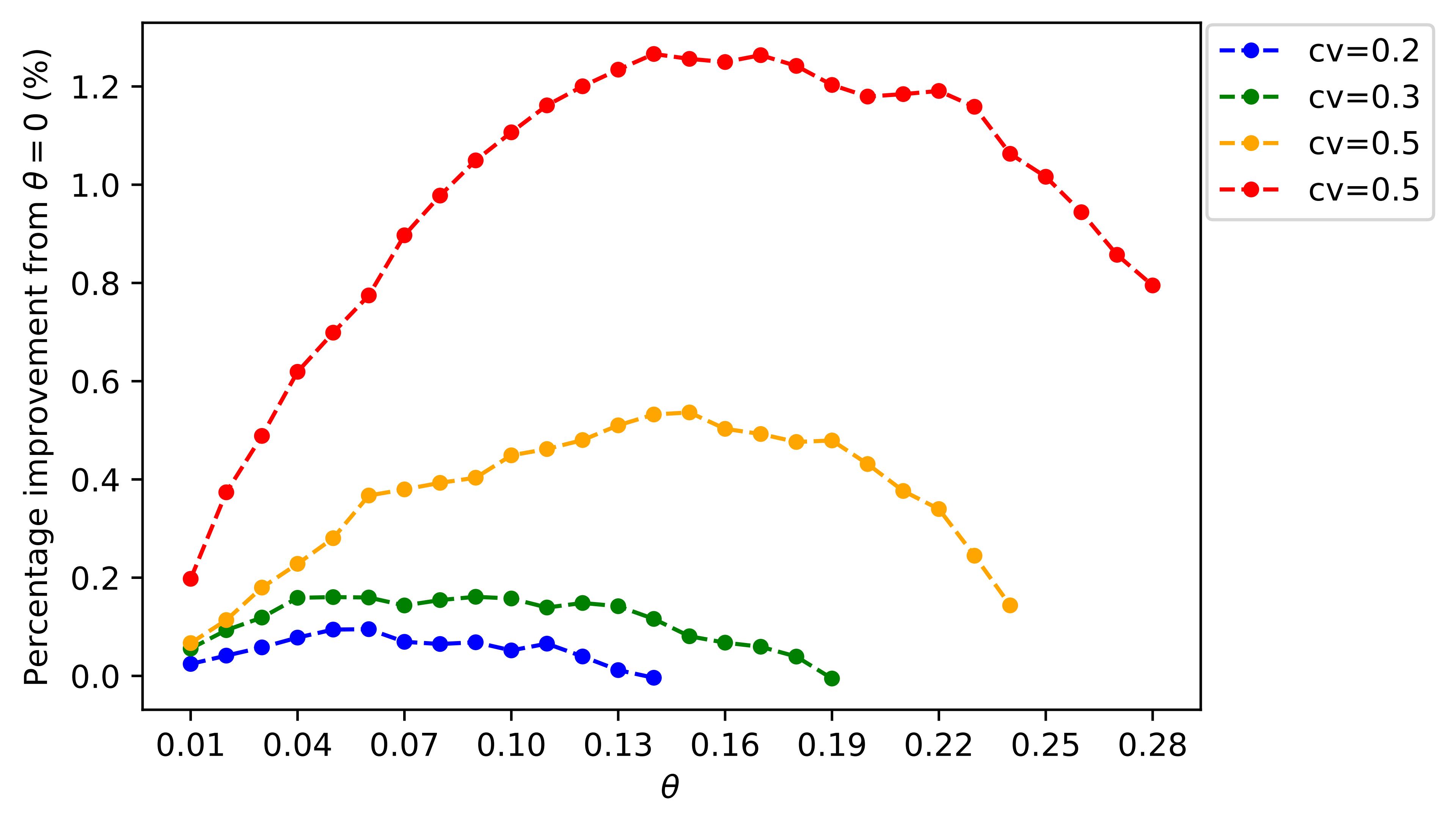}}
  \caption{Percentage improvements from WV pricing v.s. $\theta$ ($C/D=0.9$, $pf=1.25$)}
  \label{fig:cv}
\end{figure}

Finally, we discuss the impact of $cv$. We vary $cv \in \{0.2, 0.3,0.4,0.5\}$ while keeping $C/D=0.9$ and $pf=1.25$ fixed. The performance improvements over WV pricing ($\theta=0$) with respect to $\theta$ in the cases of $m=3$ and $m=27$ are shown in Figures \ref{fig:cv-m3} and \ref{fig:cv-m27}, respectively. 
The figures show that when $cv$ is small (indicating low weight and volume uncertainty), the performance improvements achieved by incorporating second-order information into the WV pricing method are relatively modest. However, as the level of uncertainty increases, the improvements become more substantial. In the case of $m=3$, the optimal value of $\theta$ increases as $cv$ increases. A similar trend is observed in the case of $m=27$, although it is relatively insignificant due to simulation error and the local stabilization of the improvement near the optimal value of $\theta$.

In general, the level of performance improvement achieved by incorporating second-order information into WV pricing increases with tighter capacity constraints, higher penalty costs, and higher uncertainty in the weight and volume of bookings. However, it also depends on specific scenarios, such as the number of booking types and the arrival process. The optimal value of $\theta$ increases as penalty cost and uncertainty increase. However, it shows little correlation with the ratio of capacity to demand. In addition, since the performance improvement remains relatively stable near the optimal value of $\theta$, finding the exact optimal value of $\theta$ may not be necessary.

\section{Conclusions}\label{sec:conclude}
This paper investigates the optimal dynamic pricing for single-leg air cargo RM during the spot sale phase. At the moment of pricing, the precise weight and volume of the cargo remain unavailable to the airline. It is only when the flight is on the brink of departure and the bookings arrive that the airline discerns the precise weight and volume. We utilize an MDP framework to model the problem, and further derive a necessary condition for its optimal pricing strategy. To break the curse of dimensionality stemming from the high-dimensional state space, we propose two categories of approximation methods. One category is based on the quantity of accepted bookings. Specifically, PQ method pools weight-volume distributions of all types of bookings into a common normal distribution and then solves a one-dimensional dynamic pricing problem. We prove that PQ's pricing strategy satisfies several desirable structural properties, which are critical for the analytical validation of the model and the acceleration of computational process. AQ method adopts an agumented pricing strategy to address PQ's extreme pricing issues in scenarios with high product heterogeneity. To further enhance the performance of dynamic pricing in high product-heterogeneity scenario, the second category of approximation methods is founded on the expected weight and volume of accepted bookings. Specifically, WV method formulates a new MDP characterized by states reflecting expected weight and volume, and utilizes the derived value function to compute the approximation pricing strategy. We establish a theoretical upper bound for the optimality gap of total expected revenue under WV method. WVS method further incorporates the second-moment information to enhance performance in the scenarios of high uncertainty. We conduct a series of numerical experiments to compare the performance of PQ, AQ, WV, and WVS pricing under various scenarios. The results demonstrate that the quantity-based methods require lower computational efforts and yield great performance when the heterogeneity of the weight and volume distributions of diverse booking types is low. Moreover, compared to PQ pricing, AQ pricing exhibits higher resilience to the booking-type heterogeneity. The weight-volume-based methods are computationally more intensive, but yield superior performance especially in the scenarios of high cargo heterogeneity. The performance improvement from WVS pricing over WV pricing is relatively significant when penalty costs and uncertainty levels are high.

For future research extension, one possible direction is the incorporation of allotments into the dynamic pricing model. Since air cargo RM comprises two phases, the allocation of capacity between allotment contracts and spot sales becomes crucial. Existing studies, such as \cite{levin2012cargo} and \cite{moussawi2014optimal}, discuss this two-phase coupling problem with capacity control policies during spot sales. Furthermore, strategic customers may become familiar with the pricing strategy during spot sales and strategically make their bookings. Exploring the game-theoretical interaction between customers and airlines in this context would also be an interesting research direction.

\ACKNOWLEDGMENT{%
This work was supported by the National Natural Science Foundation of China (72371141), China.
}

\bibliographystyle{apa}
\bibliography{ref}

%
%
%
\begin{APPENDICES}
\section{Proofs}
\subsection{Proof of Theorem \ref{prop:V_x}}
To prove the theorem, we first give a lemma.
    \begin{lemma}
        If the non-negative differentiable functions $f(x)$ and $g(x)$ (where $g(x) > 0$) satisfy $[\ln f(x)]' \leq [\ln g(x)]',$ then the function $h(x) = f(x)/g(x)$ is decreasing with respect to $x$.
    \label{lem:1}
    \end{lemma}
    \noindent Based on Equation \eqref{eq:general}, 
    \begin{equation*}
        V_t(\mathbf{x})= \sum_{i=1}^m \sup_{r_i}\Big\{m_{t}^i(1-F_t^i(r_i))\left(r_iQ_i-(V_{t+1}(\mathbf{x})-V_{t+1}(\mathbf{x}+\mathbf{e}_i))\right)\Big\} + V_{t+1}(\mathbf{x}).
    \end{equation*}
    To optimize the pricing of type $i$ bookings at state $(t,\mathbf{x})$, we take the first derivative of $m_{t}^i(1-F_t^i(r_i))\left(r_iQ_i-(V_{t+1}(\mathbf{x})-V_{t+1}(\mathbf{x}+\mathbf{e}_i))\right)$ with respect to $r_i$ and set it equal to zero. This leads to the first-order condition for the optimal pricing of type $i$ bookings, which can be expressed as Equation \eqref{eq:ropt_gen},
    \begin{equation*}
    r_i = \frac{1-F_t^i(r_i)}{f_t^i(r_i)}+\frac{V_{t+1}(\mathbf{x})-V_{t+1}(\mathbf{x}+\mathbf{e}_i)}{Q_i}.
    \end{equation*}
    Let $G_t^i(r_i)=r_i-[1-F_t^i(r_i)]/f_t^i(r_i)$. If $G_t^i(r_i)$ is an increasing function with respect to $r_i$, then the above equation has a unique solution. 
    Let $(G_t^i(r_i))' \geq 0$, we can simplify it to its equivalent form as follows, 
    \begin{equation}
        \frac{2d[\ln(1-F_t^i(r_i))]}{dr_i} \leq \frac{d[\ln f_t^i(r_i)]}{dr_i}.
        \label{eq:prop1_1}
    \end{equation}
    According to Lemma \ref{lem:1}, the inequality \eqref{eq:prop1_1} is equivalent to $(1-F_t^i(r_i))^2/f_t^i(r_i)$ being a decreasing function of $r_i$. 
    Therefore, if the function $(1-F_t^i(r_i))^2/f_t^i(r_i)$ decreases with respect to $r_i$, then the above first-order condition Equation \eqref{eq:ropt_gen} has a unique solution. Assuming that $f_t^i(r_i)$ is bounded, this unique solution is the optimal pricing.

\subsection{Proof of Lemma \ref{lem:cvx}}
    At each realization, assume the values of $X_1,X_2,...,X_{n+1}$ are $x_1,x_2,...,x_{n+1}$, and the value of $S_n$ is $s_n$. Since $g(x)$ is convex, we have 
    \begin{equation}
        \begin{aligned}
        \frac{x_n}{x_n+x_{n+1}}g(s_{n+1})+\frac{x_{n+1}}{x_n+x_{n+1}}g(s_{n-1}) &\geq g(\frac{x_n}{x_n+x_{n+1}}s_{n+1}+\frac{x_{n+1}}{x_n+x_{n+1}}s_{n-1})\\
        &=g(s_n).
        \end{aligned}
    \label{eq:cvx_1}
    \end{equation}

    Then Equation \eqref{eq:cvx_1} is also satisfied for random variables, i.e.,
    \begin{equation}
        \frac{X_n}{X_n+X_{n+1}}g(S_{n+1})+\frac{X_{n+1}}{X_n+X_{n+1}}g(S_{n-1}) \geq g(S_n).
    \label{eq:cvx_2}
    \end{equation}

Take expectations on Equation \eqref{eq:cvx_2}, we have 
\begin{equation}
    E\left[\frac{X_n}{X_n+X_{n+1}}g(S_{n+1})\right]+E\left[\frac{X_{n+1}}{X_n+X_{n+1}}g(S_{n-1})\right] \geq E[g(S_n)].
\label{eq:cvx_3}
\end{equation}

On the one side, since $X_1,X_2,...,X_{n+1}$ i.i.d., $(X_n,X_{n+1})$ is independent of $S_{n-1}$. We have
\begin{equation}
    E\left[\frac{X_{n+1}}{X_n+X_{n+1}}g(S_{n-1})\right]=E\left[\frac{X_{n+1}}{X_n+X_{n+1}}\right]E[g(S_{n-1})]=\frac{1}{2}E[h(S_{n-1})].
\end{equation}
The second equality is because of symmetry.

On the other side, 
\begin{align}
    &E\left[\frac{X_n}{X_n+X_{n+1}}g(S_{n+1})\right]-\frac{1}{2}E[g(S_{n+1})] \nonumber \\
    =& E\left[\frac{X_n-X_{n+1}}{2(X_n+X_{n+1})}g(S_{n+1})\right] \nonumber \\
    =& E\left[\frac{X_{n+1}-X_{n}}{2(X_n+X_{n+1})}g(S_{n+1})\right] \nonumber \\
    =& \frac{1}{2}\left(E\left[\frac{X_n-X_{n+1}}{2(X_n+X_{n+1})}g(S_{n+1})\right]+E\left[\frac{X_{n+1}-X_{n}}{2(X_n+X_{n+1})}g(S_{n+1})\right] \right) \nonumber \\
    =& 0,
\label{eq:cvx_4}
\end{align}
where the second equality is because of symmetry.

Combine Equations \eqref{eq:cvx_3}-\eqref{eq:cvx_4}, we have 
\begin{equation*}
    \frac{1}{2}E[h(S_{n-1})]+\frac{1}{2}E[g(S_{n+1})] \geq E[g(S_n)].
\end{equation*}
This completes the proof.

\subsection{Proof of Proposition \ref{prop:inventory}}
    For $t=T$, $V_T^{PQ}(x)-V_T^{PQ}(x+1)\geq V_T^{PQ}(x-1)-V_T^{PQ}(x)$ is equivalent to 
    \begin{align}
        & E\left[h_w\left((\sum_{l=1}^{x+1}W_l^s-C_w)^+\right) \right]- E\left[h_w\left((\sum_{l=1}^{x}W_l^s-C_w)^+\right) \right] \nonumber \\
        \geq & E\left[h_w\left((\sum_{l=1}^{x}W_l^s-C_w)^+\right) \right]- E\left[h_w\left((\sum_{l=1}^{x-1}W_l^s-C_w)^+\right) \right].
    \label{eq:cvx_E}
    \end{align}
    The volume term is omitted because of symmetry.

    Note that $h_w(z)$ and $(z-C_w)^+$ are both non-negative, non-decreasing, convex functions, then $h_w((z-C_w)^+)$ is convex. Let $g(z)=h_w((z-C_w)^+)$, and regard $W_i$ as $X_i$. By Lemma \ref{lem:cvx}, the inequality \eqref{eq:cvx_E} holds. 

    Then we take a sample path argument to prove the proposition for $t\leq T-1$, i.e., $2V_t^{PQ}(x)\geq V_t^{PQ}(x-1)+V_t^{PQ}(x+1).$ This proof is inspired by the proof of Theorem $1$ from \cite{zhao2000optimal}.

    Assuming there are four sample paths, labeled as $1,2,\bar{1},\bar{2}$. At period $t$, $x-1$ bookings have been accepted on path $1$, $x+1$ bookings have been accepted on path $2$, $x$ bookings have been accepted on path $\bar{1}$ and $\bar{2}$ respectively. 
    Let $\pi_1$ and $\pi_2$ be the nonrandomized optimal pricing strategies on Q-MDP for paths $1$ and $2$. Using these strategies, the total expected revenue from path $1$ and path $2$ is $V_t^{PQ}(x-1)+V_t^{PQ}(x+1)$ respectively. For paths $\bar{1}$ and $\bar{2}$ we construct $\bar{\pi}_1, \bar{\pi}_2$, making the total revenue from them no less than $V_t^{PQ}(x-1)+V_t^{PQ}(x+1)$. Thus, by the optimality of $V_t^{PQ}(x)$, we can complete the proof.

    Now, for a realization of the demand process on path $1$, we generate three identical copies for the other three paths. Let $r^i(\tau)=(r_1^i(\tau),\dots,r_m^i(\tau))^T$ be the realization of the pricing process under $\pi_i$, $i=1,2, \ \tau=t,t+1,...,T-1$. Starting from $\tau=t$, set $\bar{r}^1(\tau)={\rm max}\{r^1(\tau),r^2(\tau)\}=(\max\{r^1_1(\tau),r^2_1(\tau)\},...,
    \max\{r^1_m(\tau),r^2_m(\tau)\})^T$ and $\bar{r}^2(\tau)={\rm min}\{r^1(\tau),r^2(\tau)\}=(\min\{r^1_1(\tau),r^2_1(\tau)\},...,
    \min\{r^1_m(\tau),r^2_m(\tau)\})^T$ until $t_0$, the first time when the total number of accepted bookings on path $1$ is exactly one more than path $\bar{1}$. After that, let $\bar{r}^i(\tau)=r^i(\tau),\ i=1,2, \ \tau=t_0+1,t_0+2,...,T-1$. Next, it will be shown that the constructed pricing strategies $\bar{\pi}_1, \bar{\pi}_2$ can make the sum of the revenues for paths $\bar{1},\bar{2}$ no less than paths $1,2$.

    In the time interval $t$ to $t_0$, the only period path $1 (2)$ and path $\bar{1} (\bar{2})$ may generate different revenues is when a booking, such as a type $i$ booking, arrives at period $\tau$ and $r_i^1(\tau) < r_i^2(\tau)$. There are three possibilities for the booking:
    \begin{enumerate}
    \item[{\rm (1)}] The reservation price of this booking is lower than $r_i^1(\tau)$, and it is lost on all the four sample paths.
    \
    \item[{\rm (2)}] The reservation price of this booking is above $r_i^2(\tau)$, and it is accepted on all the four sample paths.
    \
    \item[{\rm (3)}] The reservation price of this booking falls between $r_i^1(\tau)$ and $r_i^2(\tau)$, and it is accepted on paths $1$,$\bar{2}$ and lost on paths $2$,$\bar{1}$.
    \end{enumerate}

    In cases (1) and (2), the sum of the revenues for paths $1$ and $2$ is equal to $\bar{1}$ and $\bar{2}$, respectively, which are $0$ and $(r_i^1(\tau)+r_i^2(\tau))Q_i$. In case (3), when $t_0=\tau$, the sum of revenues for paths $1$ and $2$ is still equal to $\bar{1}$ and $\bar{2}$, respectively, which is $r_i^1(\tau)Q_i$. After $t_0$, $\pi_i$ and $\bar{\pi}_i$ ($i=1,2$) are identical, so must be their revenues. 

    We also need to verify that the sum of the penalties at period $T$ for paths $\bar{1}$ and $\bar{2}$ does not exceed the sum for paths $1$ and $2$. If $t_0$ does occur, then in periods $t, t+1,...,t_0-1$, the number of accepted bookings on path $\bar{1}$ ($\bar{2}$) is the same as on $1$ ($2$). At period $t_0$, a booking is accepted on paths $1, \bar{2}$, while it is not accepted on path $2, \bar{1}$. After $t_0$, the total number of accepted bookings on path $\bar{1}$ ($\bar{2}$) is the same as on path $1$ ($2$). Therefore, the sum of the penalties at period $T$ for paths $\bar{1}$ and $\bar{2}$ is equal to the sum for paths $1$ and $2$. In the case $t_0$ never occurs, let $y$ represent the total number of accepted bookings on path $1$. Then, $y+2$, $y+1$, and $y+1$ respectively represent the total number of accepted bookings on path $2$, $\bar{1}$, and $\bar{2}$. Since $2V_T^{PQ}(y+1)\geq V_T^{PQ}(y)+V_T^{PQ}(y+2)$, we can show that the sum of penalties for paths $\bar{1}$ and $\bar{2}$ at period $T$ does not exceed the sum for paths $1$ and $2$.

    Summarizing the above discussion, we complete the proof.

\subsection{Proof of Proposition \ref{prop:time}}
    The proposition is equivalent to 
    \begin{equation*}
        V_t^{PQ}(x)-V_{t+1}^{PQ}(x)\geq V_t^{PQ}(x+1)-V_{t+1}^{PQ}(x+1).
    \end{equation*}

    According to Equation \eqref{eq:PQ}, we have
    \begin{equation}
        V_t^{PQ}(x)-V_{t+1}^{PQ}(x) = \sum_{i=1}^m \sup_{r_i}\big\{b_t^i(r_i)(r_iQ_i-(V_{t+1}^{PQ}(x)-V_{t+1}^{PQ}(x+1)))\big\},
    \label{eq:pr_time_1}
    \end{equation}
    \begin{equation}
        V_t^{PQ}(x+1)-V_{t+1}^{PQ}(x+1) = 
        \sum_{i=1}^m \sup_{r_i}\big\{ b_t^i(r_i)(r_iQ_i-(V_{t+1}^{PQ}(x+1)-V_{t+1}^{PQ}(x+2)))\big\}.
    \label{eq:pr_time_2}
    \end{equation}

    By Proposition \ref{prop:inventory}, we have
    \begin{equation}
        -(V_{t+1}^{PQ}(x)-V_{t+1}^{PQ}(x+1)) \geq -(V_{t+1}^{PQ}(x+1)-V_{t+1}^{PQ}(x+2)).
    \label{eq:pr_time_3}
    \end{equation}

    Since $b_t^i(r_i)\geq 0$, combine Equations \eqref{eq:pr_time_1}-\eqref{eq:pr_time_3}, we have
    \begin{equation*}
        V_t^{PQ}(x)-V_{t+1}^{PQ}(x)\geq V_t^{PQ}(x+1)-V_{t+1}^{PQ}(x+1).
    \end{equation*}

\subsection{Proof of Proposition \ref{prop:Qprice_x}}
    Define 
    \begin{align}
        h_t^i(r_i,x)&=b_t^i(r_i)(r_iQ_i-(V_{t+1}^{PQ}(x)-V_{t+1}^{PQ}(x+1))), \label{eq:h} \\
        h_t(\mathbf{r},x)&=\sum_{i=1}^m h_t^i(r_i,x)+V_{t+1}^{PQ}(x). \label{eq:hb}
    \end{align}
    
    Then by Equation \eqref{eq:PQ}, $V_t^{PQ}(x)= \sum_{i=1}^m \sup_{r_i}h_t^i(r_i,x) + V_{t+1}^{PQ}(x)$.

    We only need to prove $r_i^{PQ}(t,x)\leq r_i^{PQ}(t,x+1)$ for all $i=1,2,...,m$.

    For the pricing $\mathbf{r}(t,x+1)=(r_1(t,x+1),r_2(t,x+1),...,r_m(t,x+1))^T$, if there exists an $i\in\{1,2,...,m\}$ such that $r_i(t,x+1)<r_i^{PQ}(t,x)$, set $\mathbf{\Tilde{r}}(t,x+1)=\mathbf{r}(t,x+1)-r_i(t,x+1)\mathbf{e}_i+r_i^{PQ}(t,x)\mathbf{e}_i.$ Note that 
    \begin{align*}
        & h_t(\mathbf{\Tilde{r}}(t,x+1),x)-h_t(\mathbf{r}(t,x+1),x) \nonumber \\
        &= \sum_{j=1}^m\left[ h_t^j(\Tilde{r}_j(t,x+1),x)-h_t^j(r_j(t,x+1),x) \right] \nonumber \\
        &= h_t^i(r_i^{PQ}(t,x),x) - h_t^i(r_i(t,x+1),x) \nonumber \\
        & \geq 0,
    \end{align*}
    where the last inequality holds for the optimality of $r_i^{PQ}(t,x)$ in the Q-MDP. 

    To prove $r_i^{PQ}(t,x)\leq r_i^{PQ}(t,x+1)$, we only need to prove
    \begin{equation}
    \begin{aligned}
        h_t(\mathbf{\Tilde{r}}(t,x+1),x+1)-h_t(\mathbf{r}(t,x+1),x+1) \geq 
        h_t(\mathbf{\Tilde{r}}(t,x+1),x)-h_t(\mathbf{r}(t,x+1),x) \ (\geq 0).
    \end{aligned}
    \label{eq:rQ_1}
    \end{equation}

    Equation \eqref{eq:rQ_1} indicates for all $r_i(t,x+1)<r_i^{PQ}(t,x)$, we can adjust the pricing $r_i(t,x+1)$ to $r_i^{PQ}(t,x)$, which would increase the expected revenue-to-go at state $x+1$.

    Substitute the expression of $h_t(\mathbf{r},x)$ in Equations \eqref{eq:h}-\eqref{eq:hb} back into Equation \eqref{eq:rQ_1} and simplify, we only need to prove
    \begin{equation*}
    \begin{aligned}
         \relax [(V_{t+1}^{PQ}(x)-V_{t+1}^{PQ}(x+1))-(V_{t+1}^{PQ}(x+1)-V_{t+1}^{PQ}(x+2))]\times 
        [F_t^i(r_i(t,x+1))-F_t^i(r_i^{PQ}(t,x))] \geq 0.
    \end{aligned}
    \end{equation*}

    Since $r_i(t,x+1)<r_i^{PQ}(t,x)$, this is equivalent to $V_{t+1}^{PQ}(x+1)-V_{t+1}^{PQ}(x+2) \geq V_{t+1}^{PQ}(x)-V_{t+1}^{PQ}(x+1)$, which holds by Proposition \ref{prop:inventory}.

\subsection{Proof of Proposition \ref{prop:Qprice_t}}
    The proof is similar to the proof of Proposition \ref{prop:Qprice_x}. We only need to prove $r_i^{PQ}(t+1,x)\leq r_i^{PQ}(t,x)$ for all $i=1,2,...,m$.

    For the pricing $\mathbf{r}(t+1,x)=(r_1(t+1,x),r_2(t+1,x),...,r_m(t+1,x))^T$, if there exits an $i\in\{1,2,...,m\}$ such that $r_i(t+1,x)>r_i^{PQ}(t,x)$, set $\mathbf{\Tilde{r}}(t+1,x)=\mathbf{r}(t+1,x)-r_i(t+1,x)\mathbf{e}_i+r_i^{PQ}(t,x)\mathbf{e}_i.$ Note that
    \begin{align*}
        & h_t(\mathbf{\Tilde{r}}(t+1,x),x)-h_t(\mathbf{r}(t+1,x),x) \nonumber \\
        &= \sum_{j=1}^m\left[ h_t^j(\Tilde{r}_j(t+1,x),x)-h_t^j(r_j(t+1,x),x) \right] \nonumber \\
        &= h_t^i(r_i^{PQ}(t,x),x) - h_t^i(r_i(t+1,x),x) \nonumber \\
        & \geq 0.
    \end{align*}
    where the last inequality holds for the optimality of $r_i^{PQ}(t,x)$ in the Q-MDP.

    To prove $r_i^{PQ}(t+1,x)\leq r_i^{PQ}(t,x)$, we only need to prove
    \begin{equation}
    \begin{aligned}
        h_{t+1}(\mathbf{\Tilde{r}}(t+1,x),x)-h_{t+1}(\mathbf{r}(t+1,x),x)  \geq 
         h_t(\mathbf{\Tilde{r}}(t+1,x),x)-h_t(\mathbf{r}(t+1,x),x) \ (\geq 0).
    \label{eq:rQ_2}
    \end{aligned}
    \end{equation}

    Substitute the expression of $h_t(\mathbf{r},x)$ in Equations \eqref{eq:h}-\eqref{eq:hb} back into Equation \eqref{eq:rQ_2} and simplify, we only need to prove
    \begin{equation*}
    \begin{aligned}
         \relax [(V_{t+1}^{PQ}(x)-V_{t+1}^{PQ}(x+1))-(V_{t+2}^Q(x)-V_{t+2}^Q(x+1))]\times 
        [F^i(r_i(t+1,x))-F^i(r_i^{PQ}(t,x))] \geq 0.
    \end{aligned}
    \end{equation*}

    Since $r_i(t+1,x)>r_i^{PQ}(t,x)$, this is equivalent to $V_{t+1}^{PQ}(x)-V_{t+1}^{PQ}(x+1)\geq V_{t+2}^Q(x)-V_{t+2}^Q(x+1)$, which holds by Proposition \ref{prop:time}.

\subsection{Proof of Proposition \ref{prop:compare_1}}
    We prove the proposition by backward induction.

    First, for $t=T$, since $h_w(z),h_v(z),z^+$ are all non-decreasing convex functions with respect to $z$, $h_w(z^+),h_v(z^+)$ are non-decreasing convex functions with respect to $z$ by the properties of composition of convex functions. From Jensen's Inequality, we have
    
    \begin{align*}
        E\left[h_w\left((\sum_{i=1}^{m}\sum_{j=1}^{x_i}W_{ij}-C_w)^+\right)\right] &\geq h_w\left((\sum_{i=1}^{m}\sum_{j=1}^{x_i}E[W_{ij}]-C_w)^+\right), \\
        E\left[h_v\left((\sum_{i=1}^{m}\sum_{j=1}^{x_i}V_{ij}-C_v)^+\right)\right] &\geq h_v\left((\sum_{i=1}^{m}\sum_{j=1}^{x_i}E[V_{ij}]-C_v)^+\right).
    \end{align*}

    Adding the two inequalities, we have
    \begin{equation*}
        \bar{V}_T(\mathbf{x}) \geq V_T(\mathbf{x}).
    \end{equation*}

    Suppose the proposition is true for $t+1 \ (t \leq T-1)$, then $\bar{V}_{t+1}(\mathbf{x}) \geq V_{t+1}(\mathbf{x})$. Then

    \begin{equation*}
    \begin{aligned}
        \bar{V}_{t}(\mathbf{x}) =&\ {\rm sup_{r_i}}\Bigg\{\sum_{i=1}^mb_t^i(r_i)\left(r_iQ_i+\bar{V}_{t+1}(\mathbf{x+e_i})\right) +\left(1-\sum_{i=1}^mb_t^i(r_i)\right)\bar{V}_{t+1}(\mathbf{x}) \Bigg\} \\
        \geq & \ {\rm sup_{r_i}}\Bigg\{\sum_{i=1}^mb_t^i(r_i)\left(r_iQ_i+V_{t+1}(\mathbf{x+e_i})\right) +\left(1-\sum_{i=1}^m b_t^i(r_i)\right)V_{t+1}(\mathbf{x}) \Bigg\} \\
        = & \ V_t(\mathbf{x}),
    \end{aligned}
    \end{equation*} 
    where the second inequality is obtained using the case for $t+1$.

\subsection{Proof of Proposition \ref{prop:idealwv}}
    We prove the proposition by backward induction.
    First, for $t=T$, from Equations \eqref{eq:wv}-\eqref{eq:wv_boundary},\eqref{eq:def_wv} and \eqref{eq:CE_boundary}, we have
    
    \begin{equation*}
    \begin{aligned}
        V_T^{WV}(\mathbf{x})&=\Tilde{V}_T(\sum_{i=1}^mx_iw_i,\sum_{i=1}^mx_iv_i) \\
        &= -h_w\left((\sum_{i=1}^mx_iw_i-C_w)^+\right)-h_v\left((\sum_{i=1}^mx_iv_i-C_v)^+\right) \\
        &= \Bar{V}_T(\mathbf{x}).
    \end{aligned}
    \end{equation*}

    Suppose the proposition is true for $t+1 \ (t \leq T-1)$, then $V_{t+1}^{WV}(\mathbf{x})=\Bar{V}_{t+1}(\mathbf{x})$. Then from Equations \eqref{eq:wv}, \eqref{eq:def_wv}, we have
    
    \begin{align*}
        V_t^{WV}(\mathbf{x})=&\Tilde{V}_t(\sum_{i=1}^mx_iw_i,\sum_{i=1}^mx_iv_i) \nonumber \\
        =&\sum_{i=1}^m{\rm sup_{r_i}}\Big\{b_t^i(r_i)\Big(r_iQ_i-\big(\Tilde{V}_{t+1}(\sum_{j=1}^mx_jw_j,\sum_{j=1}^mx_jv_j) - \nonumber \\
        & \Tilde{V}_{t+1}(\sum_{j=1}^mx_jw_j+w_i,\sum_{j=1}^mx_jv_j+v_i)\big)\Big)\Big\} + \Tilde{V}_{t+1}(\sum_{j=1}^mx_jw_j,\sum_{j=1}^mx_jv_j) 
     \nonumber \\
        =&\sum_{i=1}^m{\rm sup_{r_i}}\Big\{b_t^i(r_i)\left(r_iQ_i-(V_{t+1}^{WV}(\mathbf{x})-V_{t+1}^{WV}(\mathbf{x}+\mathbf{e_i}))\right)\Big\} + V_{t+1}^{WV}(\mathbf{x}) \nonumber \\
        =& \sum_{i=1}^m{\rm sup_{r_i}}\Big\{b_t^i(r_i)\left(r_iQ_i-(\Bar{V}_{t+1}(\mathbf{x})-\Bar{V}_{t+1}(\mathbf{x}+\mathbf{e_i}))\right)\Big\} + \Bar{V}_{t+1}(\mathbf{x}) \nonumber \\
        =& \Bar{V}_t(\mathbf{x}),
    \end{align*}
    where the second from the bottom equality is obtained using the case for $t+1$.

\subsection{Proof of Theorem \ref{prop:DCE}}
    The proof of this theorem will employ the following lemma.
    \begin{lemma}[\cite{Gallego1994}]
    For any random variable $X$ with finite mean $\mu$ and finite standard deviation $\sigma$, and for any real number $a$,

    \begin{equation*}
        E[(X-a)^+]\leq \frac{\sqrt{\sigma^2+(a-\mu)^2}-(a-\mu)}{2}.
    \end{equation*}
    \label{lem:cite}
    \end{lemma}

    Let $\mathbf{r}^*(t,\mathbf{x})=(r_1^*(t,\mathbf{x}),...,r_m^*(t,\mathbf{x}))^T$ denote the optimal pricing obtained from solving the general model, and $\Bar{\mathbf{r}}(t,\mathbf{x})=(\Bar{r}_1(t,\mathbf{x}),...,\Bar{r}_m(t,\mathbf{x}))^T$ denote the CE pricing. Then, $\Bar{V}_t(\mathbf{x})$ and $J_t^{CE}(\mathbf{x})$ satisfy the following equations:

    \begin{align}
    \Bar{V}_t(\mathbf{x})=&\sum_{i=1}^mb_t^i(\bar{r}_i)\left(\Bar{r}_iQ_i-(\Bar{V}_{t+1}(\mathbf{x})-\Bar{V}_{t+1}(\mathbf{x}+\mathbf{e_i}))\right)+\Bar{V}_{t+1}(\mathbf{x}), \nonumber \\ 
    &\quad\quad\quad t=0,1,2,...,T-1.  \label{eq:VCE} \\
    \bar{V}_T(\mathbf{x})=&-h_w\left((\sum_{i=1}^{m}x_iw_i-C_w)^+\right)-h_v\left((\sum_{i=1}^{m}x_iv_i-C_v)^+\right), \label{eq:VCE_boundary} \\
    J_t^{CE}(\mathbf{x})=&\sum_{i=1}^mb_t^i(\bar{r}_i)\left(\Bar{r}_iQ_i-(J^{CE}_{t+1}(\mathbf{x})-J^{CE}_{t+1}(\mathbf{x}+\mathbf{e_i}))\right)+J^{CE}_{t+1}(\mathbf{x}), \nonumber \\ 
    &\quad\quad\quad t=0,1,2,...,T-1. \label{eq:JCE} \\
    J^{CE}_T(\mathbf{x})=&-E\left[h_w\left((\sum_{i=1}^{m}\sum_{j=1}^{x_i}W_{ij}-C_w)^+\right)+h_v\left((\sum_{i=1}^{m}\sum_{j=1}^{x_i}V_{ij}-C_v)^+\right)\right]. \ \label{eq:JCE_bound}
\end{align}

    Now we prove the theorem by backward induction.

    First, for $t=T$, using Equations \eqref{eq:VCE_boundary} and \eqref{eq:JCE_bound}, and according to Lemma \ref{lem:cite}, we have

    \begin{align*}
        D_T^{CE}(\mathbf{x})=&\Bar{V}_T(\mathbf{x})-J_T^{CE}(\mathbf{x}) \nonumber \\
        =& E\left[h_w\left(\sum_{i=1}^{m}\sum_{j=1}^{x_i}W_{ij}-C_w\right)^+\right] - h_w\left(\sum_{i=1}^{m}x_iw_i-C_w\right)^+ \nonumber \\
        &+ E\left[h_v\left(\sum_{i=1}^{m}\sum_{j=1}^{x_i}V_{ij}-C_v\right)^+\right]-h_v\left(\sum_{i=1}^{m}x_iv_i-C_v\right)^+ \nonumber \\
        \leq & \frac{h_w}{2}\left(\sqrt{\sigma_W^2(\mathbf{x})+(C_w-\Bar{w}(\mathbf{x}))^2}-|C_w-\Bar{w}(\mathbf{x})|\right) \nonumber \\
        &+ \frac{h_v}{2}\left(\sqrt{\sigma_V^2(\mathbf{x})+(C_v-\Bar{v}(\mathbf{x}))^2}-|C_w-\Bar{v}(\mathbf{x})|\right) \nonumber \\
        \leq & \max\limits_{||\mathbf{x}||_1\leq x_{max}}\{ M(\mathbf{x})\}. \label{eq:DCE_T}
    \end{align*}

    Suppose the theorem is true for $t+1 \ (t \leq T-1)$, then $D_{t+1}^{CE}(\mathbf{x})\leq \max\limits_{||\mathbf{x}||_1\leq x_{max}}\{ M(\mathbf{x})\}$.
    Then using Equations \eqref{eq:VCE} and \eqref{eq:JCE}, we have
    
    \begin{align*}
        D_t^{CE}(\mathbf{x})=&\Bar{V}_t(\mathbf{x})-J_t^{CE}(\mathbf{x}) \nonumber \\
        =&\sum_{i=1}^mb_t^i(\Bar{r}_i)(\Bar{V}_{t+1}(\mathbf{x}+\mathbf{e}_i)-J^{CE}_{t+1}(\mathbf{x}+\mathbf{e}_i)) \nonumber \\
        & + \left(1-\sum_{i=1}^mb_t^i(\Bar{r}_i) \right)(\Bar{V}_{t}(\mathbf{x})-J^{CE}_{t}(\mathbf{x})) \nonumber \\
        =& \sum_{i=0}^m b_t^i(\Bar{r}_i) D_{t+1}^{CE}(\mathbf{x}+\mathbf{e}_i)  \\
        \leq & \max\limits_{||\mathbf{x}||_1\leq x_{max}}\{ M(\mathbf{x})\},
    \end{align*}
    where $\mathbf{e}_0=\mathbf{0}_m,\ \sum_{i=0}^m b_t^i(\Bar{r}_i)=1$. The last inequality is obtained using the case for $t+1$.

\end{APPENDICES}




\end{document}